\documentclass[11pt]{article}
\usepackage{latexsym,amssymb,amsmath,amsthm,graphicx,float,bm,grffile,hyperref}

\newcommand{\pd}{\partial}

\newcommand{\R}{\mathbb{R}}

\newcommand{\<}{\langle}
\renewcommand{\>}{\rangle}

\newcommand{\tr}{\mbox{tr}}

\newcommand{\beq}{\begin{equation}}
\newcommand{\eeq}{\end{equation}}
\renewcommand{\tilde}{\widetilde}

\newcommand{\bit}{\begin{itemize}}
\newcommand{\eit}{\end{itemize}}
\newcommand{\ben}{\begin{enumerate}}
\newcommand{\een}{\end{enumerate}}
\newcommand{\x}{\mathbf{x}}
\newcommand{\y}{\mathbf{y}}
\newcommand{\bc}{\mathbf{c}}

\newtheorem{definition}{Definition}

\newtheorem{theorem}{Theorem}

\title{Compressed absorbing boundary conditions via matrix probing}
\author{Rosalie B\'elanger-Rioux and Laurent Demanet\footnote{Department of Mathematics, MIT. RBR is supported by the National Sciences and Engineering Research Council of Canada. This work was also supported by AFOSR, ONR, NSF, Total SA, and the Alfred P. Sloan Foundation.}}

\date{January 2014}

\begin{document}
\maketitle

\begin{abstract}
Absorbing layers are sometimes required to be impractically thick in order to offer an accurate approximation of an absorbing boundary condition for the Helmholtz equation in a heterogeneous medium. It is always possible to reduce an absorbing layer to an operator at the boundary by layer-stripping elimination of the exterior unknowns, but the linear algebra involved is costly. We propose to bypass the elimination procedure, and directly fit the surface-to-surface operator in compressed form from a few exterior Helmholtz solves with random Dirichlet data. The result is a concise description of the absorbing boundary condition, with a complexity that grows slowly (often, logarithmically) in the frequency parameter.
 
\end{abstract}
{ \bf Acknowledgments}

The authors would like to thank Thomas Hagstrom for useful discussions. 
\section{Introduction}

This paper investigates arbitrarily accurate realizations of absorbing (a.k.a. open, radiating) boundary conditions (ABC) for the 2D acoustic high-frequency Helmholtz equation in certain kinds of heterogeneous media. Instead of considering a specific modification of the partial differential equation, such as a perfectly matched layer, we study the broader question of compressibility of the nonlocal kernel that appears in the exact boundary integral form of the ABC. 

The full boundary integral viewpoint invites to rethink ABCs as a 2-step numerical scheme, where 
\begin{enumerate}
\item a precomputation sets up an expansion of the kernel of the boundary integral equation, then 
\item a fast algorithm is used for each application of this integral kernel at the open boundaries in a Helmholtz solver.
\end{enumerate}
This two-step approach may pay off in scenarios when the precomputation is amortized over a large number of solves of the original equation with different data.

This paper addresses the first, precomputation step: we describe a basis for the efficient expansion of the integral kernel of the ABC in some simple 2D settings, and discuss a randomized probing procedure to quickly find the coefficients in the expansion.

\subsection{Setup}\label{sec:setup}

We consider the Helmholtz equation in $\mathbb{R}^2$,
\begin{equation}\label{eq:HE}
\Delta u(\x)+\frac{\omega^2}{c^2(\x)} u(\x) = f(\x), \qquad \x = (x_1, x_2),
\end{equation}
with compactly supported $f$. Throughout the paper we consider the unique solution determined by the Sommerfeld radiation condition (SRC) at infinity: when $c(\x)$ extends to a constant $c$ outside of a bounded set, the SRC is \cite{McLean}
\begin{equation}\label{src}
\lim_{r \rightarrow \infty} r^{1/2} \left( \frac{\pd u}{\pd r} - ik u \right) = 0, \qquad k = \frac{\omega}{c},
\end{equation}
where $r$ is the radial coordinate.

We then seek to reformulate the SRC on the boundary $\pd \Omega$, so that the resulting solution inside $\Omega$ matches that of the free-space problem \eqref{eq:HE}, \eqref{src}. Let $G(\x,\y)$ be the fundamental solution for this problem, i.e., the solution when $f(\x) = \delta(\x-\y)$. Define the single and double layer potentials, respectively, on some closed contour $\Gamma$ by the following, for $\psi, \ \phi$ on $\Gamma$ (see details in \cite{McLean}, \cite{CK}):
\[
 S \psi (\x)=\int_{\Gamma} G(\x,\y) \ \psi(\y) \ dS_y, \qquad T \phi (\x) =\int_{\Gamma} \frac{\pd G}{\pd \nu_{\y}} (\x,\y) \ \phi(\y) \ dS_{\y},
\]
where $\nu$ is the outward pointing normal to the curve $\Gamma$, and $\x$ is not on $\Gamma$.
Now let $u^+$ satisfy the Helmholtz equation \eqref{eq:HE} in the exterior domain $\R \setminus \overline{\Omega}$, along with the SRC \eqref{src}. Then Green's third identity is satisfied in the exterior domain: using $\Gamma= \pd \Omega$, we get
\begin{equation}\label{eq:GRF}
T u^+ - S \frac{\pd u}{\pd \nu}^+ =  u^+, \qquad \x \in \R^2 \setminus \overline{\Omega}.
\end{equation}
Finally, using the jump condition of the the double layer $T$, we obtain Green's identity on the boundary $\pd \Omega$:
\[
(T - \frac{1}{2} I ) \, u^+ - S \frac{\pd u}{\pd \nu}^+ =  0, \qquad \x \in \pd \Omega.
\]

When the single-layer potential $S$ is invertible\footnote{This is the case when there is no interior resonance at frequency $\omega$, which could be circumvented by the use of combined field integral equations as in \cite{CK}. The existence and regularity of $D$ ultimately does not depend on the invertibility of $S$.}, we can let $D = S^{-1} (T - \frac{1}{2} I )$, and equivalently write (dropping the $+$ in the notation)
\begin{equation}\label{eq:dtn-abc}
\frac{\pd u}{\pd \nu} = D u, \qquad \x \in \pd \Omega.
\end{equation}
The operator $D$ is called the exterior Dirichlet-to-Neumann map (DtN). It is independent of the right hand side $f$ of \eqref{eq:HE} as long as $f$ is supported in $\Omega$. 
The notion that (\ref{eq:dtn-abc}) can serve as an exact ABC was made clear in a homogeneous medium, e.g., in \cite{engmaj} and in \cite{kelgiv}. Equation (\ref{eq:dtn-abc}) continues to hold even when $c(\mathbf{x})$ is heterogenous in the vicinity of $\partial \Omega$, provided the correct (often unknown) Green's function is used. The medium is indeed heterogeneous near $\partial \Omega$ in many situations of practical interest, such as in geophysics.

The numerical realization of ABC typically involves absorbing layers that become impractical for difficult $c(\mathbf{x})$, or for high accuracy. We instead propose to realize the ABC by directly compressing the integral kernel of $D$, so that the computational cost of its setup and application would become competitive when (\ref{eq:HE}) is to be solved multiple times. Hence this paper is not concerned with the design of a new ABC, but rather with the reformulation of existing ABCs that otherwise require a lot of computational work per solve. In many situations of practical interest we show that it is possible to ``learn" the integral form of $D$, as a precomputation, from a small number of solves of the exterior problem with the expensive ABC. By ``small number", we mean a quantity essentially independent of the number of discretization points $N$ -- in practice as small as 1 or as large as 50. We call this strategy matrix probing.

For much of the paper, the letter $D$ refers to the exact discrete realization of the DtN map, while $\tilde{D}$ is used for its approximation in the probing framework. We begin with a brief overview of absorbing conditions and why they may become computationally unwieldy in heterogeneous media. Next, we review the complexity of existing methods for forming the matrix $D$ and solving (\ref{eq:HE}) with \eqref{eq:dtn-abc}. We then present matrix probing as a compression strategy for $D$, and show numerical results that document its complexity. We also show results on the solution of the Helmholtz equation using $\tilde{D}$ for $D$ in \eqref{eq:dtn-abc}. Finally, we prove convergence of our algorithm for approximating $D$ in the special case of the half plane in a uniform medium.

Note that the interior DtN map has a much higher degree of complexity than the exterior one, because it needs to encode all the waves that travel along the broken geodesic rays that go from one part of the domain $\Omega$ to another (i.e., rays bouncing inside the domain.) In contrast, the exterior DtN map rarely needs to take into account multiple scattering if the solution is outgoing. We only consider the exterior DtN map in this work, and refer to it as the DtN map for simplicity.

\subsection{Discrete absorbing boundary conditions}\label{sec:related}

There are many ways to realize an absorbing boundary condition for the wave or Helmholtz equation. Some ABCs are surface-to-surface, such as in  \cite{engmaj}, \cite{kelgiv}, \cite{higdon}, \cite{hag}. Others involve surrounding the computational domain $\Omega$ by an absorbing layer (\cite{berenger},  \cite{druskin}, \cite{appcol}). This approach is desirable because the parameters of the layer can usually be adjusted to obtain a desired accuracy. 

While a layer should preferably be as thin, its design involves at least two different factors: 1) waves that enter the layer must be significantly damped before they re-enter the computationaly domain, and 2) reflections created when waves cross the domain-layer interface must be minimized. The Perfectly Matched Layer of B\'erenger (called PML, see \cite{berenger}) is a convincing solution to this problem in a uniform acoustic medium. Its performance often carries through in a general heterogeneous acoustic medium $c(\mathbf{x})$, though its derivation strictly speaking does not.
We may still define a layer-based scheme from a transformation of the spatial derivatives which mimics the one done in a homogeneous medium, by replacing the Laplacian operator $\Delta$ by some $\Delta_{layer}$ inside the PML, but this layer will not be perfectly matched anymore. In this case, reflections from the interface between $\Omega$  and the layer are usually not small. In fact, the layer might even cause the solution to grow exponentially inside it, instead of forcing it to decay (\cite{diazjoly}, \cite{back}). It has been shown in \cite{adiabatic} that, in some cases of interest to the optics community with nonuniform media, PML for Maxwell's equations can still work, but the layer needs to be made very thick in order to minimize reflections at the interface. In this case, the Helmholtz equation has to be solved in a very large computational domain, where most of the work will consist in solving for the PML. In the setting where \eqref{eq:HE} has to be solved a large number of times, a precomputation to speed up the application of the PML (or any other accurate but slow ABC) might be of interest.

Discrete absorbing layers may need to be quite wide in practice. Call $L$ this width (in meters). Although this is not a limitation of the framework presented in this paper, we discretize the Helmholtz operator in the most elementary way using the standard five-point difference stencil. Put $h = 1/N$ for the grid spacing, where $N$ Is the number of points per dimension for the interior problem, inside the unit square $\Omega = [0,1]^2$. 

While $\Omega$ contains $N^2$ points, the total number of unknowns is $O\left((N+2w)^2\right)$ in the presence of the layer, where $w=L/h$ is its width in number of grid points. In a uniform medium, the PML width $L$ needed is a fraction of the wavelength, i.e. $L \sim \lambda=\frac{2\pi}{\omega} \sim \frac{1}{N}$, so that we need a constant number of points independently of $N$: $w=L/h=LN \sim 1$. However, in nonuniform media, the heterogeneity of $c(\mathbf{x})$ can limit the accuracy of the layer. If we consider an otherwise uniform medium with an embedded scatterer outside of $\Omega$, then the PML will have to be large enough to enclose this scatterer, no matter $N$. For more general, heterogeneous media such as the ones considered in this paper, we often observe that convergence as a function of $L$ or $w$ is delayed compared to a uniform medium. That means that we have $L \sim L_0$ so that $w \sim NL_0$ or $w = O(N)$, as we assume in the sequel.

In the case of a second-order discretization, the rate at which one must increase $N$ in order to preserve a constant accuracy in the solution, as $\omega$ grows, is about $N \sim \omega^{1.5}$. This unfortunate phenomenon, called the pollution effect, is well-known: it begs to \emph{increase} the resolution, or number of points per wavelength, of the scheme as $\omega$ grows  \cite{nvsom,BabPollut}. As we saw, the width of the PML may be as wide as a constant value $L_0$ independent of $N$, hence its width generally needs to scale as $O(\omega^{1.5})$ grid points.

Before we explain our approach for compressing an ABC, we explain the most straightforward way of obtaining a DtN map from an ABC, by eliminating the unknowns in the absorbing layer in order to obtain a reduced system on the interior nodes. This solution, however, is computationally impractical.

\subsection{Layer-stripping for the Dirichlet-to-Neumann map}\label{sec:strip}

We write the system for the discrete Helmholtz equation as
\begin{equation}\label{HEsys}
\begin{pmatrix} \\
 & A & P & \\ 
 & & & \\
& P^T& C & \\
& & &
\end{pmatrix}  \quad 
\begin{pmatrix} \\ u_{out} \\ \\ u_{\Omega} \\ \\ \end{pmatrix}
= \begin{pmatrix} \\ 0 \\ \\ f_\Omega \\ \\ \end{pmatrix},
\end{equation}
with $A=\Delta_{layer} + k^2 I$ and $C=\Delta + k^2 I$, with $\Delta$ overloaded to denote discretization of the Laplacian operator, and $\Delta_{layer}$ the discretization of the Laplacian operator inside the PML layer. We wish to eliminate the exterior unknowns $u_{out}$ from this system in order to have a new system which only depends on the interior unknowns $u_{\Omega}$. The most obvious way of eliminating those unknowns is to form the Schur complement $S=C-P^TA^{-1}P$ of $A$ by any kind of Gaussian elimination. For instance, in the standard raster scan ordering of the unknowns, the computational cost of this method is $O(w^4)$ --- owing from the fact that $A$ is a sparse banded matrix of size $O(w^2)$ and band at least $2w$. Alternatively, elimination of the unknowns can be performed by layer-stripping, starting with the outermost unknowns from $u_{out}$, until we eliminate the layer of points that is just outside of $\pd \Omega$. The computational cost will be $O(w^4)$ in this case as well. To see this, let $u_{w}$ be the points on the outermost layer, $u_{w-1}$ the points in the layer just inside of $u_{w}$, etc. Then we have the following system:
\[
\begin{pmatrix} \\
 & A_w & P_w & \\ 
 & & & \\
& P_w^T& C_w & \\
& & & \\
\end{pmatrix}  \quad 
\begin{pmatrix} \\ u_{w} \\ \\ \vdots \\ \\ \end{pmatrix}
= \begin{pmatrix} \\ 0 \\ \\ \vdots  \\ \\ \end{pmatrix}
\]

Note that, because of the five-point stencil, $P_w$ has non-zeros exactly on the columns corresponding to $u_{w-1}$. Hence the matrix $P_w^TA_w^{-1}P_w$ in the first Schur complement $S_w=C_w-P_w^TA_w^{-1}P_w$ is non-zero exactly at the entries corresponding to $u_{w-1}$. It is then clear that, in the next Schur complement, to eliminate the next layer of points, the matrix $A_{w-1}$ (the block of $S_w$ corresponding to the points $u_{w-1}$) to be inverted will be full. For the same reason, every matrix $A_j$ to be inverted thereafter, for every subsequent layer to be eliminated, will be a full matrix. Hence at every step the cost of forming the corresponding Schur complement is at least on the order of $m^3$, where $m$ is the number of points in that layer. Hence the total cost of eliminating the exterior unknowns by layer stripping is approximately
\[ \sum_{j=1}^{w} (4(n+2j))^3 = O(w^4). \]

Similar arguments can be used for the Helmholtz equation in 3 dimensions. In this case, the computational complexity of the Schur complement or layer-stripping methods would be $O(w^3 (w^2)^2)=O(w^7)$ or $~\sum_{j=1}^{w} (6(n+2j)^2)^3=O(w^7)$. Therefore, direct elimination of the exterior unknowns is quite costly. Some new insight is needed to construct the DtN map more efficiently.

We now remark that, whether we eliminate exterior unknowns in one pass or by layer-stripping, we obtain a reduced system. It looks just like the original Helmholtz system on the interior unknowns $u_{\Omega}$, except for the top left block, corresponding to $u_{0}$ the unknowns on $\pd \Omega$, which has been modified by the elimination procedure. Hence with the help of some dense matrix $D$ we may write the reduced, $N^2$ by $N^2$ system as
\begin{equation}\label{HEred}
Lu=
\begin{pmatrix} 
(hD-I)/h^2  & I/h^2 & 0 & & \cdots & \\ 
 &  &  & \\
 I/h^2 &  & & \\
       & & & \\
0 & & [ \; \Delta + k^2 I \; ] & \\
& & & \\
\vdots & & & \\
\\
\end{pmatrix} \quad 
\begin{pmatrix} u_{0} \\  \\   u_{-1}  \\ \\   u_{-2}  \\ \\ \vdots \\ \\ \end{pmatrix} 
= \begin{pmatrix}  0  \\  \\ f_{-1} \\ \\ f_{-2} \\ \\ \vdots \\ \\ \end{pmatrix} 
\end{equation}
and we have thus obtained an absorbing boundary condition which we may use on the boundary of $\Omega$, independent of the right-hand-side $f$. Indeed, if we call $u_{-1}$ the first layer of points inside $\Omega$, we have $ (I-hD)u_{0} = u_{-1} $, or
\[
 \frac{u_{0} - u_{-1}}{h}=Du_{0} ,
\]
a numerical realization of the DtN map in \eqref{eq:dtn-abc}, using the ABC of choice, say PML. Indeed, elimination can be used to reformulate any computationally intensive ABC, not just absorbing layers, into a realization of \eqref{eq:dtn-abc}. Any ABC is equivalent to a set of equations relating unknowns on the surface to unknowns close to the surface, and possibly auxiliary variables. Again, elimination can reduce those equations to relations involving only unknowns on the boundary and on the first layer inside the boundary, to obtain a numerical DtN map $D$.

The reduced system (\ref{HEred}) is smaller than the original system (\ref{HEsys}) and often faster to solve. A drawback is that forming this matrix $D$ by elimination is prohibitive, as we have just seen. Instead, this paper
suggests adapting the framework of matrix probing in order to obtain $D$ in reasonable complexity.

\subsection{Matrix probing for the Dirichlet-to-Neumann map}\label{sec:probing-intro}

The idea of matrix probing is that a matrix $D$ with adequate structure can sometimes be recovered from the knowledge of a fixed, small number of matrix-vector products $Dg_k$, where $g_k$ are typically random vectors. In the case where $D$ is the DtN map, each $g_k$ consists of Dirichlet data on $\pd \Omega$, and each application $Dg_k$ requires solving an exterior Helmholtz problem. 

The dimensionality of $D$ needs to be limited for recovery from a few $Dg_k$ to be possible, but matrix probing is \emph{not} an all-purpose low-rank approximation technique. Instead, it is the property that $D$ has an efficient representation in some adequate pre-set basis that makes recovery from probing possible. As opposed to the randomized SVD method which requires the number of matrix-vector applications to be greater than the rank \cite{Halko-randomSVD}, matrix probing can recover interesting structured operators from a single matrix-vector application \cite{Chiu-probing, Demanet-probing}.

In order to describe the structure of the DtN map $D$, notice first that $D$ has a 4 by 4 block structure if, for example, $\pd \Omega$ has 4 sides. Hence $D$ is $n \times n$ where $n=4N$.  As an integral kernel, $D$ would have singularities at the junctions between these blocks (due to the singularities in $\pd \Omega)$, so we respect this feature by probing $D$ block by block.

We now describe a model for $M$, any $N \times N$ block of $D$, that will sufficiently lower its dimensionality to make probing possible. Assume we can write $M$ as
\begin{equation}\label{eq:Dexp}
M \approx \sum_{j=1}^p c_j B_j
\end{equation}
where the $B_j$'s are fixed, known basis matrices, that need to be chosen carefully in order to give an accurate approximation of $M$. In case when the medium $c$ is homogeneous, we typically let $B_j$ be a discretization of the integral kernel
\begin{equation}\label{eq:Bj}
B_j(x,y)= \frac{e^{ik|x-y|}}{(h+|x-y|)^{j/2}},
\end{equation}
where $h=1/N$ is the discretization parameter. We usually add another index to the $B_j$, and a corresponding multiplicative factor, to allow for a smooth dependence on $x+y$ as well. 
We shall further detail our choices and discuss their rationales in Section \ref{sec:probe}. 

Given a random vector $z \sim N(0,I_N)$ (other choices are possible), the product $w=Mz$ and the expansion \eqref{eq:Dexp}, we can now write
\begin{equation}\label{mp}
w=Mz \approx \sum_{j=1}^p c_j B_j z = \Psi_z \, \bc.
\end{equation}
Multiplying this equation on the left by the pseudo-inverse of the $N$ by $p$ matrix $\Psi_z$ will give an approximation to $\bc$, the coefficient vector for the expansion \eqref{eq:Dexp} of $M$. More generally, if several applications $w_k = M z_k$, $k = 1,\ldots, q$ are available, a larger system is formed by concatenating the $\Psi_{z_k}$ into a tall-and-thin $Nq$ by $p$ matrix ${\bm \Psi}$. The computational work is dominated, here and in other cases \cite{Chiu-probing, Demanet-probing}, by the matrix-vector products $Dg$, or $Mz$.

In a nutshell, recovery of $\bc$ works under mild assumptions on $B_j$, and when $p$ is a small fraction of $Nq$ up to log factors. More details about the performance guarantees of probing are given in Section \ref{sec:probe}.

The advantage of the specific choice of basis matrix (\ref{eq:Bj}), and its generalizations explained in the sequel, is that it results in accurate expansions with a number of parameters $p$ which is ``essentially independent" of $N$, namely that ot grows either logarthmically in $N$, or at most like a very sublinear fractional power law (such as $N^{0.12}$.) This is in sharp contract to the scaling for the layer width, $w = O(N)$ grid points, discussed earlier.

Let us now argue theoretically why $p$ indeed depends very weaky on $N$, in the special case of the half-space DtN map in a uniform medium.

\subsection{Rate of convergence of the basis matrices expansion to the half-space DtN map}

The form of $B_j$ suggested in equation (\ref{eq:Bj}) is motivated by the fact that they provide a good expansion basis for the uniform medium half-space DtN map in $\R^2$. For the Helmholtz equation $\Delta u(x,y) +k^2 u(x,y) = 0$ in $y > 0$ with the Sommerfeld radiation condition and $c\equiv 1$ (thus in this section and in Section \ref{sec:pf}, $\omega=k$), we have
\begin{equation}\label{dtn}
 \left. \partial_y u(x,y) \right|_{y=0} = \int_{-\infty}^\infty K(|x-x'|) u(x',0) \ dx',
\end{equation}
where
\begin{equation}\label{kernel}
 K(r)=\frac{ik}{2r} H_1^{(1)}(kr),
\end{equation}
with $H_1^{(1)}$ the Hankel function of the first kind, of order 1. We let $D$ for the operator mapping $u(\cdot,0)$ to $\left. \partial_y u(x,y) \right|_{y=0}$.

Since $K(r)$ is singular at $r = 0$, and since discretization effects dominate near the diagonal in the matrix representation of the DtN map, we only study the representation of $K$ in the range $r_0 \leq r \leq 1$, with $r_0$ on the order of $1/k$. Let $\tilde{K}(r) = K \chi_{[ \frac{1}{k}, 1 ]}(r)$. Denote by $\tilde{D}$ the corresponding operator with integral kernel $\tilde{K}(|x-x'|)$. Since the diagonal is cut out at level $r_0$, we also modify the basis in a cosmetic but convenient way by replacing $\frac{e^{ikr}}{(r_0+r)^{j/\alpha}}$ by $\frac{e^{ikr}}{r^{j/\alpha}}$.

\begin{theorem}\label{teo:main}
Let $\alpha > \frac{2}{3}$, and let $\tilde{K}_p(r)$ be the best uniform approximation of $\tilde{K}(r)$ in 
\[
\mbox{span} \{ \frac{e^{ikr}}{r^{j/\alpha}} : j = 1, \ldots, p, \mbox{ and } r_0 \leq r \leq 1 \}.
\]
Assume that $r_0 = C/k$ for some $C>0$ independent of $k$. Denote by $\tilde{D}_p$ the operator defined with $\tilde{K}_p$ in place of $\tilde{K}$. Then, in the operator norm,
\[
\| \tilde{D} - \tilde{D}_p \| \leq C_\alpha \, p^{1 - \lfloor 3\alpha/2 \rfloor} \, \| \tilde{K} \|_\infty,
\]
for some $C_\alpha > 0$ depending on $\alpha$.
\end{theorem}

The proof, and a numerical illustration, are in Section \ref{sec:pf}. Growing $\alpha$ does not automatically result in a better approximation error, because a careful analysis of the proof would show that $C_\alpha$ grows factorially with $\alpha$. This behavior translates into a slower onset of convergence in $p$ when $\alpha$ is taken large, as the numerics show, which can in turn be interpreted as the result of ``overcrowding" of the basis by very look-alike functions.

Notice that $D$ is not bounded in $L^2$, but $\tilde{D}$ is after the diagonal is cut out. It is easy to see that the operator norm of $\tilde{D}$ grows like $k$, for instance by applying $\tilde{D}$ to the function $e^{-ikx}$. The uniform norm of $\tilde{K}$, however, grows like $k^2$, so the result above shows that we incur an additional factor $k$ in the error (somewhat akin to numerical pollution) in addition to the factor $k$ that we would have gotten from $\| \tilde{D} \|$.

The important point of the theorem is that the quality of approximation is otherwise independent of $k$, i.e., the number $p$ of basis functions does not need to grow like $k$ for the error to be small. In other words, it is unnecessary to ``mesh at the wavelength level" to spell out the degrees of freedom that go in the representation of the DtN map's kernel.

The goal of the paper is to approximate $D$ in more general cases than the half-space, but the result in Theorem \ref{teo:main} points the way for the design of basis matrices in the general case in Section \ref{sec:basis}.

\section{Algorithms}

\subsection{Setup for the exterior problem}\label{sec:ext}

The exterior problem is the heterogeneous-medium Helmholtz equation at angular frequency $\omega$, outside $\Omega=[0,1]^2$, with Dirichlet boundary condition $u_0=g$ on $\pd \Omega$. As in the introduction, this problem is solved numerically with the five-point stencil of finite differences. The PML starts at a fixed, small distance away from $\Omega$, so that we keep a small strip around $\Omega$ where the equations are unchanged. Recall that the width of the PML is in general as large as $O(\omega^{1.5})$ grid points. We number the sides of $\pd \Omega$ counter-clockwise starting from $(0,0)$, hence side 1 is the bottom side $(0,y)$, $0\leq y \leq 1$, side 2 is the right side $(x,1)$, $0\leq x \leq 1$, etc.

The method by which the system for the exterior problem is solved is immaterial in the scope of this paper, though for reference, the experiments in this paper use UMFPACK's sparse direct solver \cite{UMFPACK}. For treating large problems, a better solver should be used, such as the sweeping preconditioner of Engquist and Ying \cite{Hsweep,Msweep}, the shifted Laplacian preconditioner of Erlangga \cite{erlangga}, the domain decomposition method of Stolk \cite{stolk}, or the direct solver with spectral collocation of Martinsson, Gillman and Barnett \cite{dirfirst,dirstab}. This in itself is a subject of ongoing research which we shall not discuss further. 

For a given $g$, we solve the system and obtain a solution $u$ in the exterior computational domain. In particular we consider $u_{1}$, the solution in the layer just outside of $\pd \Omega$. We know from Section \ref{sec:strip} that $u_1$ and $g$ are related by $ u_{1}=(I+hD)g$ or
\begin{equation}\label{eq:D}
 \frac{u_{1} - g}{h}=Dg
\end{equation}
The matrix $D$ that this relation defines needs not be interpreted as a first-order approximation of the continous DtN map: it is the algebraic object of interest that will be ``probed" from repeated applications to different vectors $g$.

Similarily, for probing the $(i_M,j_M)$ block of $D$ -- that we generically call $M$ -- one needs matrix-vector products of $D$ with vectors $g$ of the form $[z, 0, 0, 0]^T$, $[0, z, 0, 0]^T$, etc., to indicate that the Dirichlet BC is $z$ on the side indexed by $j_M$, and zero on the other sides. The application $Dg$ is then restricted to side $i_M$.

\subsection{Matrix probing}\label{sec:probe}

We saw how to compute the matrix product $Mz$ in the previous section, where $M$ is a block of the discrete DtN map $D$. Recall we have assumed for matrix probing that we can write $M$ using $p$ known basis matrices $B_j$ as $M \approx \sum_{j=1}^p c_j B_j$. This lead to the approximation $w=Mz \approx \sum_{j=1}^p c_j B_jz = \Psi_z \bc$. Note that both $\Psi_z$ and the resulting coefficient vector $\bc$ depend on the vector $z$. In the sequel we let $z$ be gaussian iid random.

In order to improve the conditioning in taking the pseudo-inverse of the matrix $\Psi_z$ and reduce the error in the coefficient vector $\bc$, one may use $q > 1$ random realizations of $M$, that is, $w^1=Mz^1, \ldots, w^q=Mz^q$. Then, $w$ will be a long column vector containing the concatenation of the $w^i$'s, and $\Psi_z$ will have size $Nq$ by $p$, and we still solve for $\bc$ in $w=\Psi_z \bc$. There is a limit to the range of $p$ for which this system is well-posed: past work by one of us \cite{Chiu-probing} covers the precise conditions on $p$, $n$, and the following two parameters, called ``weak condition numbers'', for which recoverability of $\bc$ is accurate with high probability.

\begin{definition}
\[ \lambda = \max_j \frac{\| B_j \|_2 \sqrt{N}}{\| B_j \|_F} \]
\end{definition}
\begin{definition}\label{kap}
\[ \kappa = \mbox{cond}( N), \ N_{j \ell} = \mbox{Tr} \, (B_j^T B_\ell)\]
\end{definition}

It is desirable to have a small $\lambda$, which translates into a high rank condition on the basis matrices, and a small $\kappa$, which translates into a Riesz basis condition on the basis matrices. Having small weak condition numbers will guarantee a small failure probabilty of matrix probing and a bound on the condition number of $\Psi_z$, i.e. guaranteed accuracy in solving for $\bc$. Also, using $q > 1$ allows to use a larger $p$, to achieve greater accuracy. These results are contained in the following theorem.
\begin{theorem} (Chiu-Demanet, \cite{Chiu-probing}) Let $z$ be a Gaussian i.i.d. random vector of length $qN$, and $\Psi_z$ as above. Then $\mbox{cond}(\Psi_z) \leq 2\kappa + 1$ with high probability provided that $p$ is not too large, namely 
\[
q N \geq C \, p \, (\kappa \lambda \log N)^2,
\]
for some number $C > 0$.
\end{theorem}

As noted previously, the work necessary for probing the matrix $M$ is on the order of $q$ solves of the original problem. Indeed, computing $Mz^1, \ldots , Mz^q$ means solving $q$ times the exterior problem with the PML. This is roughly equivalent to solving the original Helmholtz problem with the PML $q$ times, assuming $w$ is at least as large as $N$. Then, computing the $qp$ products of the $p$ basis matrices with the $q$ random vectors amounts to a total of at most $qpN^2$ work, or less if the basis matrices have a fast matrix-vector product. And finally, computing the pseudo-inverse of $\Psi_z$ has cost $Nqp^2$. Hence, as long as $p,q \ll N$, the dominant cost of matrix probing comes from solving $q$ times the exterior problem with a random Dirichlet boundary condition. In our experiments, $q=O(1)$ and $p$ can be as large as $1000$ for high accuracy.

Finally, we note that the information from the $q$ solves can be re-used for any other block which is in the same block column as $M$. However, if it is needed to probe blocks of $D$ which are not all in the same block column, then another $q$ solves need to be performed, with a Dirichlet boundary condition on the appropriate side of $\partial \Omega$. This of course increases the total number of solves. Another option would be to probe all of $D$ at once, using a combination of basis matrices that have the same size as $D$, but that are 0 except on the support of each distinct block in turn. In this case, $\kappa$ remains the same because we still orthogonalize our basis matrices, but $\lambda$ doubles ($\| B_j \|_2 $ and $\| B_j \|_F$ do not change but $N \rightarrow 4N$) and this makes the conditioning worse, in particular a higher value of $q$ is needed for the same accuracy, given by $p$. Hence we have decided not to investigate further this approach, which might become more advantageous in the case of a more complicated polygonal domain.

\subsection{Choice of basis matrices}\label{sec:basis}

The essential information of the DtN map needs to be summarized in broad strokes in the basis matrices $B_j$, with the details of the numerical fit left to the probing procedure. In the case of $D$, most of its physics is contained in its \emph{oscillations} and \emph{diagonal singularity}, as predicted by geometrical optics.

A heuristic argument to obtain the form of $D$ starts from the Green's formula (\ref{eq:GRF}), that we differentiate one more time in the normal direction. After accounting for the correct jump condition, we get an alternative Steklov-Poincare identity, namely
\[
D = (T^* + \frac{1}{2} I)^{-1} H,
\]
where $H$ is the hypersingular integral operator with kernel $\frac{\pd^2 G}{\pd \nu_{\x} \pd \nu_{\y}}$, where $G(\x,\y)$ is the Green's function and $\nu_{\x}$, $\nu_{\y}$ are the normals to $\pd \Omega$ in $\x$ and $\y$ respectively. The presence of $(T^* + \frac{1}{2} I)^{-1}$ is somewhat inconsequential to the form of $D$, as it involves solving a well-posed second-kind integral equation. As a result, the properties of $D$ are qualitatively similar to those of $H$. (The exact construction of $D$ from $G$ is of course already known in a few special cases, such as the uniform-medium half-space problem considered earlier.)

In turn, geometrical optics reveals the form of $G$. In a context where there is no multi-pathing, that is, where there is a single traveltime $\tau(\x,\y)$ between any two points $\x,\y \in \Omega$, one may write a high-$\omega$ asymptotic series for $G$ as
\begin{equation}\label{eq:geoopts}
 G(\x,\y) \sim e^{i\omega \tau(\x,\y)} \sum_{j\geq 0} A_j(\x,\y) \omega^{-j},
\end{equation}
$\tau(\x,\y)$ is the traveltime between points $\x$ and $\y$, found by solving the eikonal equation
\begin{equation} \label{eq:tau}
 \| \nabla_{\x} \tau(\x,\y) \| = \frac{1}{c(\x)},
\end{equation}
and the amplitudes $A_j$ satisfy transport equations. In the case of multi-pathing (possible multiple traveltimes between any two points), the representation \eqref{eq:geoopts} of $G$ becomes instead
 \[ 
 G(\x,\y) \sim \sum_j e^{ i \omega \tau_j(\x,\y)} \sum_{k \geq 0} A_{jk}(\x,\y) \omega^{-k}, 
 \]
where the $\tau_j$'s are the traveltimes, each obeying \eqref{eq:tau} away from caustic curves. The amplitudes are singular at caustic curves in addition to the diagonal $\x=\y$, and contain the information of the Maslov indices. Note that traveltimes are symmetric: $\tau_j(\x,\y)=\tau_j(\y,\x)$, and so is the kernel\footnote{The proof of the symmetry of $D$ was shown in a slightly different setting here \cite{symm} and can be adapted to our situation.} of $D$.

The singularity of the amplitude factor in \eqref{eq:geoopts}, at $\x = \y$, is $O \left( \log | \x - \y| \right)$ in 2D and $O \left( | \x - \y |^{-1} \right)$ in 3D. After differentiating twice to obtain $H$, the homogeneity on the diagonal becomes $O \left( | \x - \y|^{-2} \right)$ in 2D and $O \left( | \x - \y |^{-3} \right)$ in 3D. For the decay at infinity, the scalings are different and can be obtained from Fourier analysis of square root singularities; the kernel of $H$ decays like $O \left(| \x - \y|^{-3/2} \right)$ in 2D, and $O \left(| \x - \y|^{-5/2} \right)$ in 3D. In between, the amplitude is smooth as long as the traveltime is single-valued.

Much more is known about DtN maps, such as the many boundedness and coercivity theorems between adequate fractional Sobolev spaces (mostly in free space, with various smoothness assumptions on the boundary). We did not attempt to leverage these fine properties of $D$ in the scheme presented here.

For all these reasons, we define the basis matrices $B_j$ as follows. Assume $\tau$ is single-valued. In 1D, denote the tangential component of $\x$ by $x$, and similarly that of $\y$ by $y$, in coordinates local to each edge with $0 \leq x,y \leq 1$. Each block $M$ of $D$ relates to a couple of edges of the square domain. Let $j = (j_1, j_2)$ with $j_1, j_2$ nonnegative integers. The general forms that we consider are
\[
\beta_j(x,y) = e^{i \omega \tau(x,y)} (h + |x-y|)^{-\frac{j_1}{\alpha}} (h + \theta(x,y))^{-\frac{j_2}{\alpha}}
\]
and
\[
\beta_j(x,y) = e^{i \omega \tau(x,y)} (h + |x-y|)^{-\frac{j_1}{\alpha}} (h + \theta(x,y))^{j_2},
\]
where $h$ is the grid spacing of the FD scheme, and $\theta(x,y)$ is an adequate function of $x$ and $y$ that depends on the particular block of interest. The more favorable choices for $\theta$ are those that respect the singularities created at the vertices of the square; we typically let $\theta(x,y) = \min(x+y, 2-x-y)$. The parameter $\alpha$ can be taken to be equal to 2, a good choice in view of the numerics and in the light of the asymptotic behaviors on the diagonal and at infinity discussed earlier.

If several traveltimes are needed for geometrical reasons, then different sets of $\beta_j$ are defined for each traveltime. (More about this in the next section.) The $B_j$ are then obtained from the $\beta_j$ by QR factorization within each block, where orthogonality is defined in the sense of the Frobenius inner product $\< A, B \> = \tr(A B^T)$. This automatically sets the $\kappa$ number of probing to 1.

In many of our test cases it appears that the ``triangular" condition $j_1 + 2 j_2 < $ constant works well. The number of couples $(j_1,j_2)$ satisfying this relation will be $p/T$, where $p$ is the number of basis matrices in the matrix probing algorithm and $T$ is the number of distinct traveltimes. The eventual ordering of the basis matrices $B_j$ respects the increase of $j_1 + 2 j_2$.

\subsection{Traveltimes}\label{sec:tt}

Determining the traveltime(s) $\tau(\x,\y)$ is the more ``supervised" part of this method, but is needed to keep the number $p$ of parameters small in the probing expansion. A few different scenarios can arise.

\begin{itemize}
\item In the case when $\nabla c(\x)$ is perpendicular to a straight segment of the boundary, locally, then this segment is itself a ray and the waves can be labeled as interfacial, or ``creeping". The direct traveltime between any two points $\x$ and $\y$ on this segment is then simply given by the line integral of $1/c(\x)$. An infinite sequence of additional interfacial waves result from successive reflections at the endpoints of the segment, with traveltimes predicted as follows.

We still consider the exterior problem for $[0,1]^2$. We are interested in the traveltimes between points $\x, \y$ on the same side of $\pd \Omega$ -- for illustration, let $\x=(x,0)$ and $\y=(y,0)$ on the bottom side of $\Omega=[0,1]^2$, with $x \leq y$ (this is sufficient since traveltimes are symmetric). Assume that all the waves are interfacial. The first traveltime $\tau_1$ corresponds to the direct path from $\x$ to $\y$. The second arrival time $\tau_2$ will be the minimum traveltime corresponding to: either starting at $\x$, going left, reflecting off of the $(0,0)$ corner, and coming back along the bottom side of $\pd \Omega$, past $\x$ to finally reach $\y$; or starting at $\x$, going past $\y$, reflecting off of the $(1,0)$ and coming straight back to $\y$. The third arrival time $\tau_3$ is the maximum of those two choices. The fourth arrival time then corresponds to starting at $\x$, going left, reflecting off of the $(0,0)$ corner, travelling all the way to the $(1,0)$ corner, and then back to $\y$. The fifth arrival time corresponds to leaving $\x$, going to the $(1,0)$ corner this time, then back to the $(0,0)$ corner, then on to $\y$. And so on. To recap, we have the following formulas:
\begin{eqnarray*}
\tau_1(\x,\y)&=& \int_x^y \frac{1}{c(t,0)} \ dt, \\
\tau_2(\x,\y)&=& \tau_1(\x,\y) + 2\min \left( \int_0^x \frac{1}{c(t,0)} \ dt, \int_y^1 \frac{1}{c(t,0)} \ dt \right), \\
\tau_3(\x,\y)&=& \tau_1(\x,\y) + 2\max \left( \int_0^x \frac{1}{c(t,0)} \ dt, \int_y^1 \frac{1}{c(t,0)} \ dt \right) = 2\int_0^1 \frac{1}{c(t,0)} \ dt - \tau_2(\x,\y), \\
\tau_4(\x,\y)&=& 2\int_0^1 \frac{1}{c(t,0)} \ dt - \tau_1(\x,\y), \\
\tau_5(\x,\y)&=& 2\int_0^1 \frac{1}{c(t,0)} \ dt + \tau_1(\x,\y), \qquad \mbox{etc.} \\
\end{eqnarray*}
All first five traveltimes can be expressed as a sum of $\pm \tau_1$, $\pm \tau_2$ and the constant phase $2\int_0^1 \frac{1}{c(t,0)} \ dt$, which does not depend on $\x$ or $\y$. In fact, one can see that any subsequent traveltime corresponding to traveling solely along the bottom boundary of $\pd \Omega$ should be again a combination of those quantities. This means that if we use $\pm \tau_1$ and $\pm \tau_2$ in our basis matrices, we are capturing all the traveltimes relative to a single side, which helps to obtain higher accuracy for probing the diagonal blocks of $D$.

This simple analysis can be adapted to deal with creeping waves that start on one side of the square and terminate on another side, which is important for the nondiagonal blocks of $D$.

\item In the case when $c(\x)$ increases outward in a smooth fashion, we are also often in presence of body waves, going off into the exterior and coming back to $\pd \Omega$. The traveltime for these waves needs to be solved either by a Lagrangian method (solving the ODE for the rays), or by an Eulerian method (solving the eikonal PDE shown earlier). In this paper we used the fast marching method of Sethian \cite{sethart} to deal with these waves in the case that we label ``slow disk" in the next section. 

\item In the case when $c(\x)$ has singularities in the exterior domain, each additional reflection creates a traveltime that should (ideally) be predicted. Such is the case of the ``diagonal fault" example introduced in the next section, where a straight jump discontinuity of $c(\x)$ intersects $\pd \Omega$ at a non-normal angle:  we can construct by hand the traveltime corresponding to a path leaving the boundary at $\x$, reflecting off of the discontinuity and coming back to the boundary at $\y$.  More precisely, we consider again $\x=(x,0)$, $\y=(y,0)$ and $x \leq y$, with $x$ larger than or equal to the $x$ coordinate of the point where the reflector intersects the bottom side of $\pd \Omega$. We then  reflect the point $\y$ across the discontinuity into the new point $\y'$, and calculate the Euclidean distance between $\x$ and $\y'$. To obtain the traveltime, we then divide this distance by the value $c(\x)=c(\y)$ of $c$ on the right side of the discontinuity, assuming that value is constant. This body traveltime is used in the case of the ``diagonal fault", replacing the quantity $\tau_2$ that was described above. This increased accuracy by an order of magnitude, as mentioned in the numerical results of the next section.

\end{itemize}

\subsection{Solving the Helmholtz equation with a compressed ABC}

Once we have obtained approximations $\tilde{M}$ of each block $M$ in compressed form through the coefficients $\bc$ using matrix probing, we construct block by block the approximation $\tilde{D}$ of $D$ and use it in a solver for the Helmholtz equation on the domain $\Omega=[0,1]^2$, with the boundary condition
$$\frac{\pd u}{\pd \nu}  = \tilde{D}u , \qquad x \in \pd \Omega.$$

\section{Numerical experiments}

\begin{figure}[H]
\begin{minipage}[b]{0.32\linewidth}
\includegraphics[scale=.32]{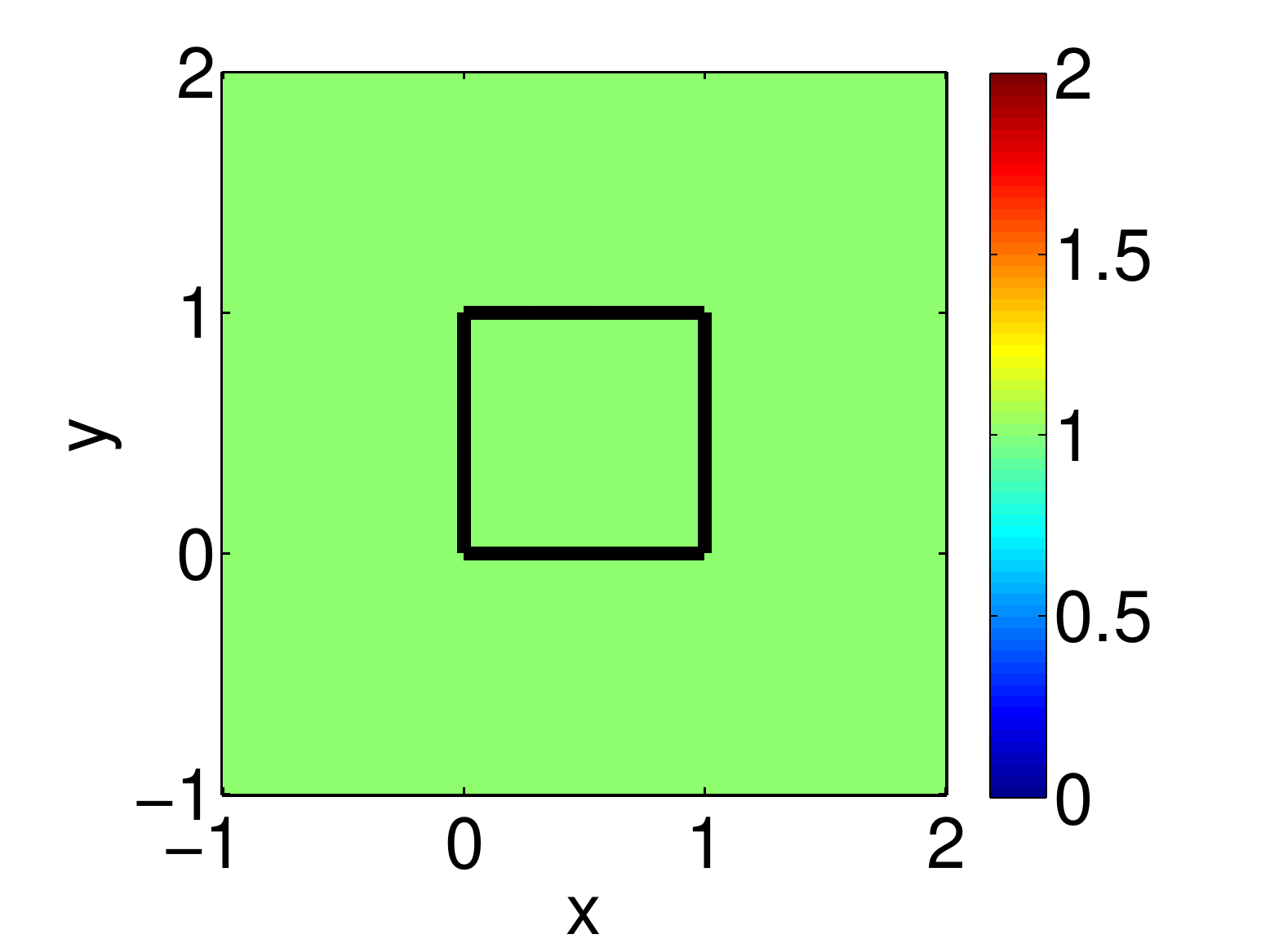}
\caption{Uniform medium.}\label{c1}
\includegraphics[scale=.32]{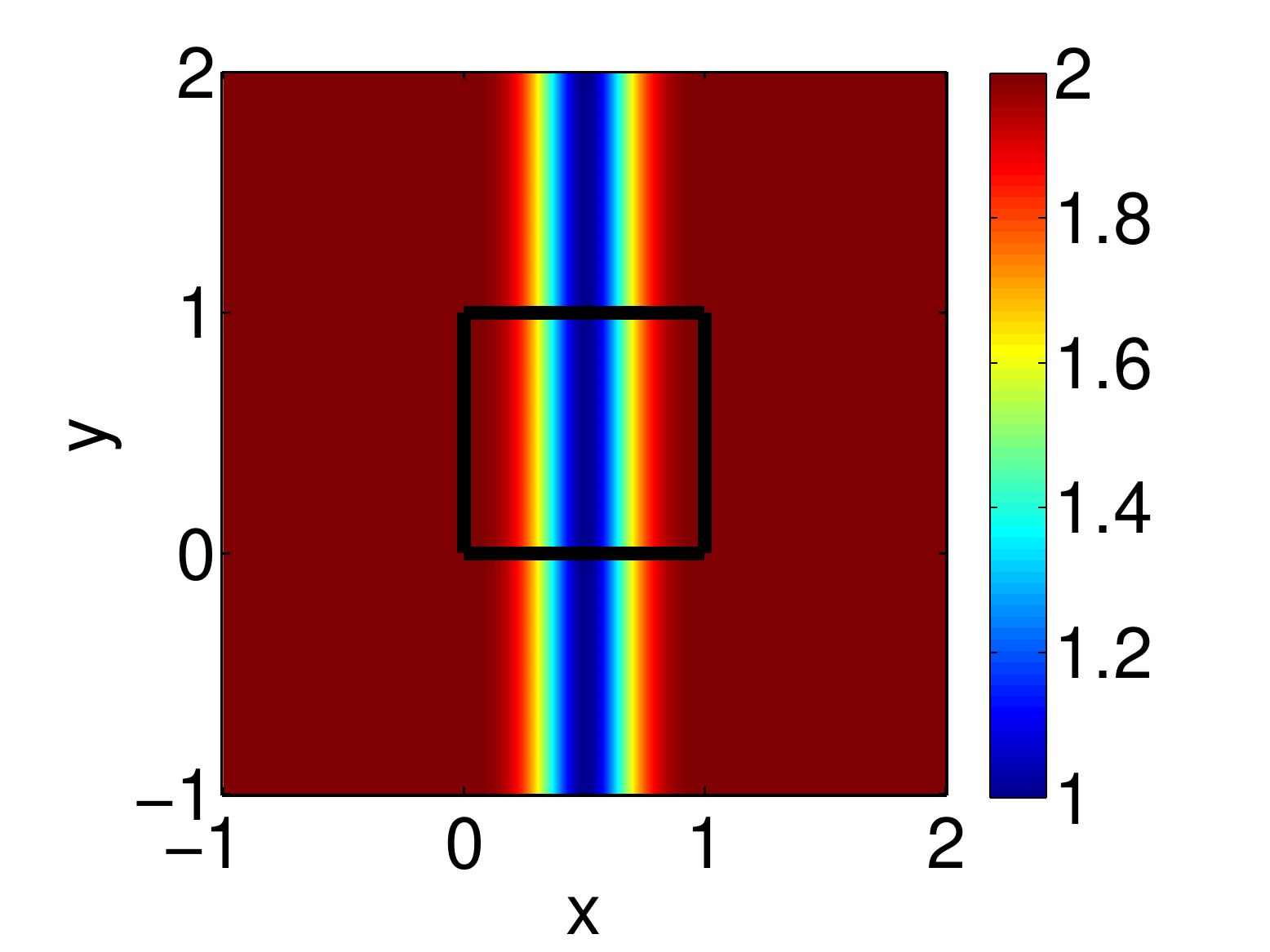}
\caption{Gaussian waveguide.}\label{wg}
\end{minipage}
\begin{minipage}[b]{0.32\linewidth}
\includegraphics[scale=.32]{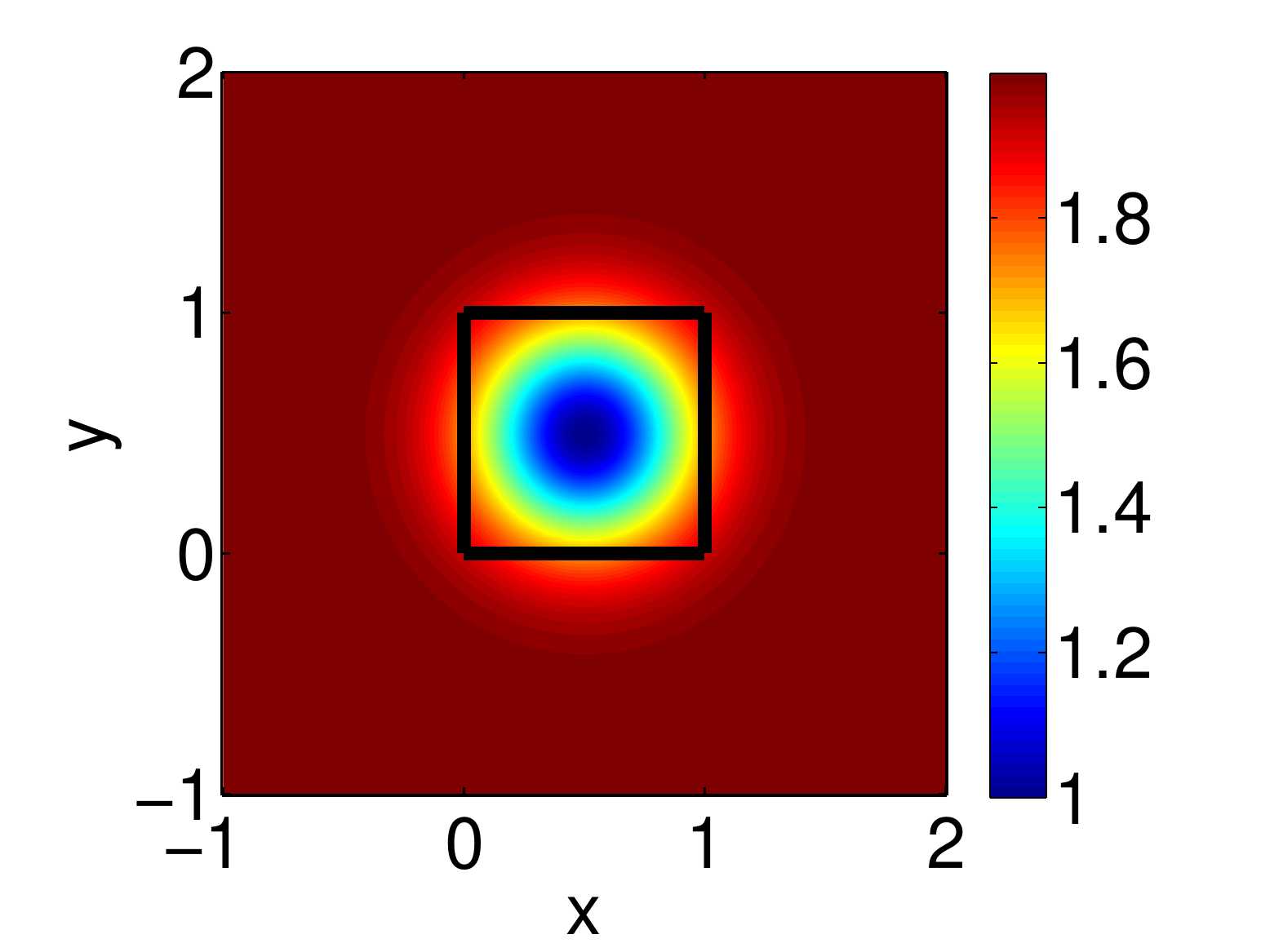}
\caption{Gaussian slow disk.}\label{slow}
\includegraphics[scale=.32]{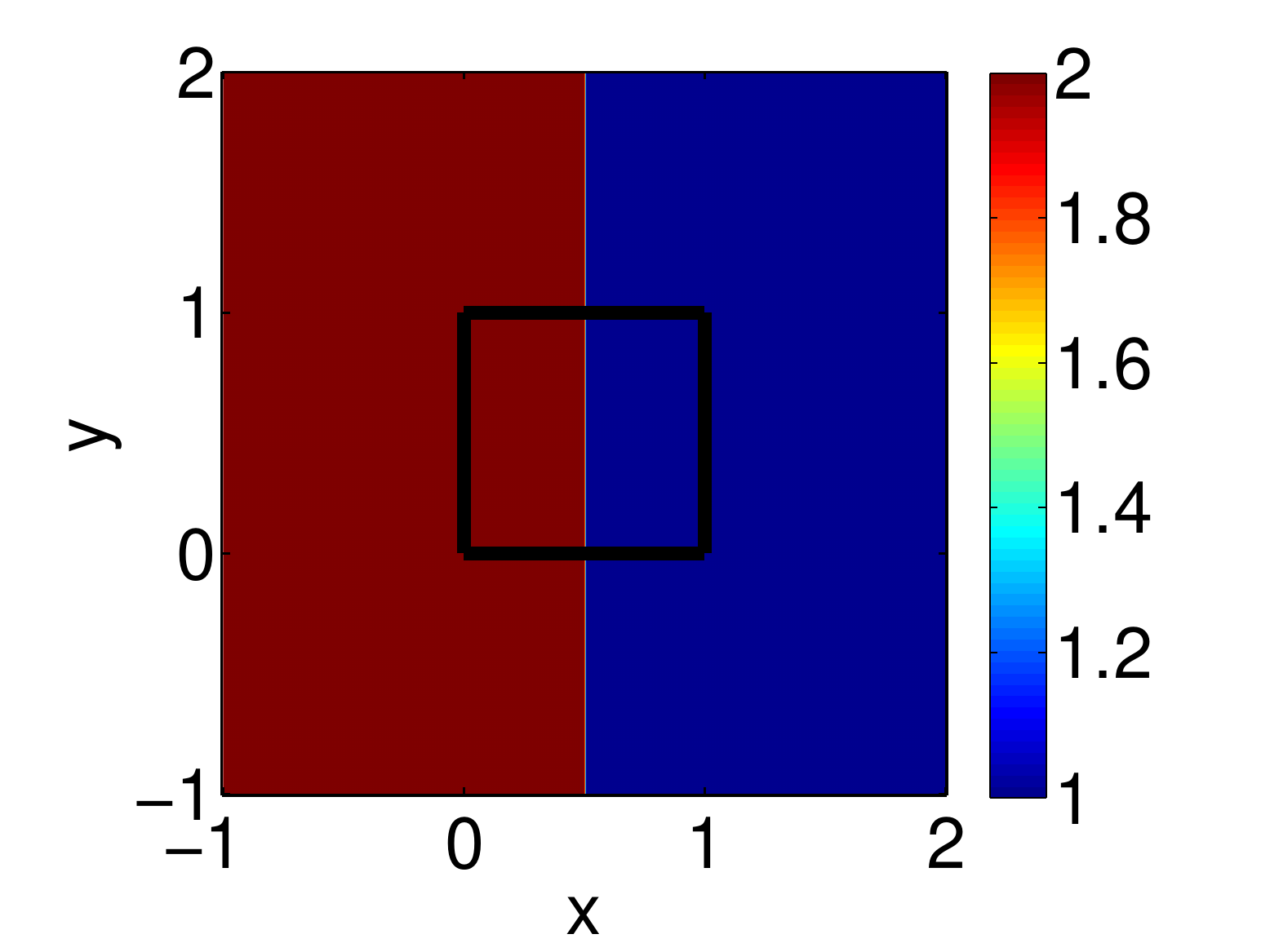}
\caption{Vertical fault.}\label{fault}
\end{minipage}
\begin{minipage}[b]{0.32\linewidth}
\includegraphics[scale=.32]{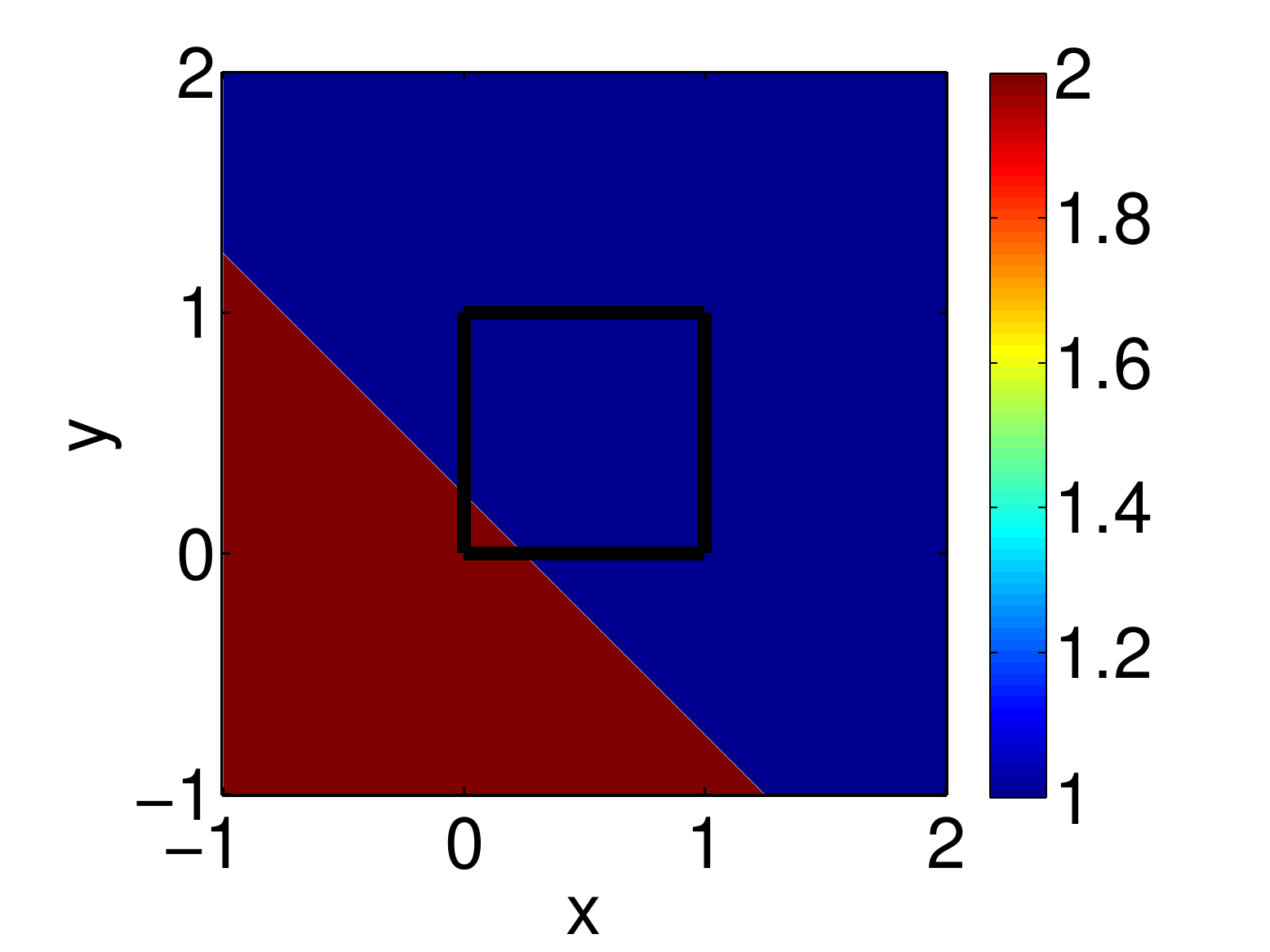}
\caption{Diagonal fault.}\label{diagfault}
\includegraphics[scale=.32]{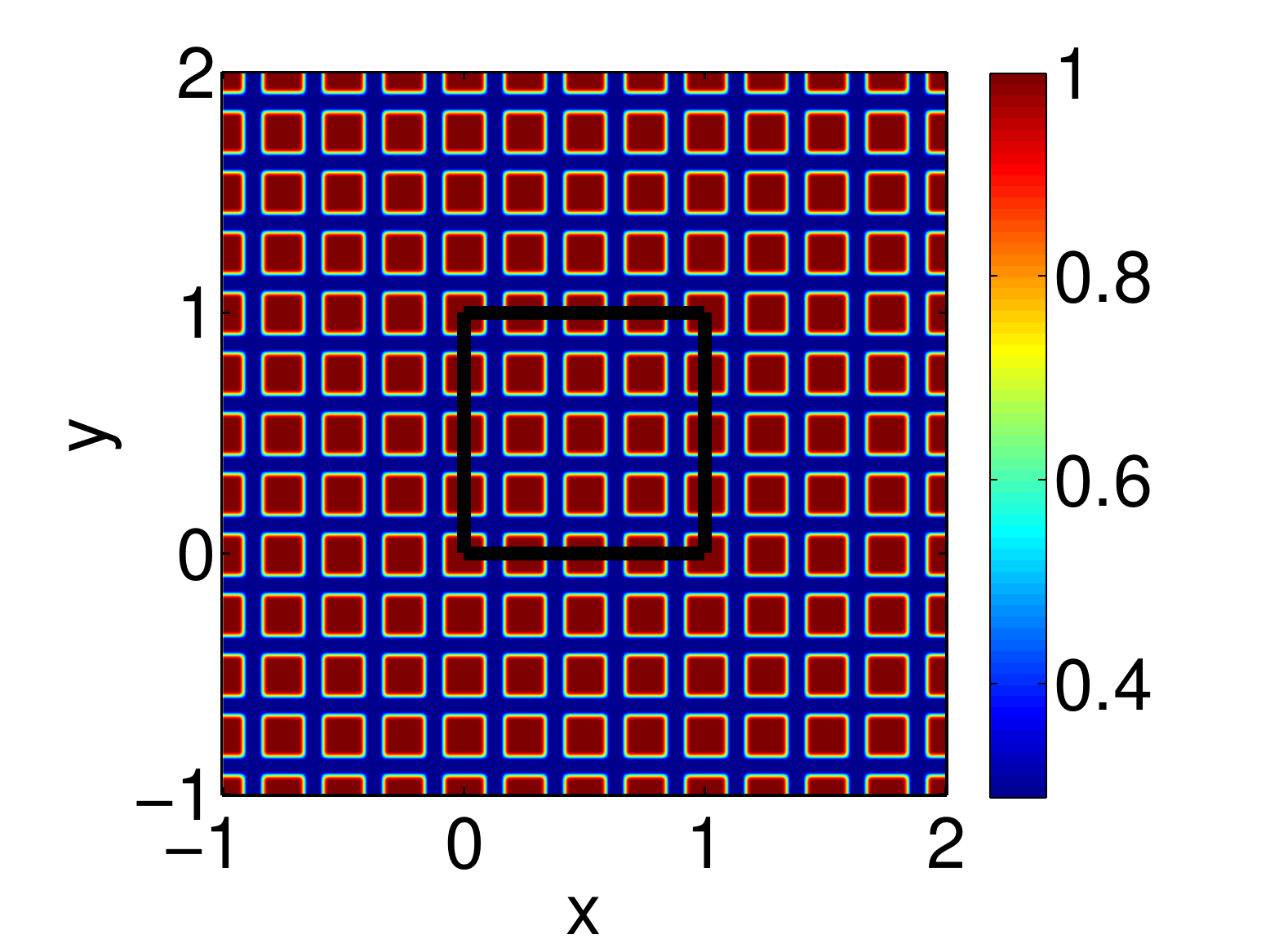}
\caption{Periodic medium.}\label{period}
\end{minipage}
\end{figure}

Our benchmark media $c(\x)$ are as follows:

\begin{enumerate}
\item a constant wave speed of 1 (Figure \ref{c1}),
\item a ``Gaussian waveguide" (Figure \ref{wg}),
\item a ``Gaussian slow disk" (Figure \ref{slow}) large enough to encompass $\Omega$ - this will cause some waves going out of $\Omega$ to come back in, 
\item a ``vertical fault" (Figure \ref{fault}), 
\item a ``diagonal fault" (Figure \ref{diagfault}), 
\item and a periodic medium (Figure \ref{period}). The periodic medium consists of square holes of velocity 1 in a background of velocity $1/\sqrt{12}$. 
\end{enumerate}
All media used are continued in the obvious way (i.e., they are \emph{not} put to a homogeneous constant) outside of the domain in which they are shown in the figures if needed. The outline of the $[0,1]^2$ box is shown in black.

We can use a standard Helmholtz equation solver to estimate the relative error in the Helmholtz equation solution caused by the Finite Difference discretization, and also the error caused by using the specified PML width. Those errors are presented in Table \ref{FDPMLerr}, along with the main parameters used in the remaining of this section, including the position of the point source or right-hand side $f$.

\begin{table}
\begin{center} \footnotesize
\begin{tabular}{|l|l|l|l|l|l|l|} \hline  
Medium			&$N$	&$\omega/2\pi$	&FD error	&$w$		& $P$	&Source position	\\ \hline  
$c \equiv 1$		&1023	&51.2		&{$2.5e-01$}	&{$4$}		&8	&$(0.5,0.25)$		\\ \hline  
waveguide		&1023	&51.2		&{$2.0e-01$}	&{$4$}		&56	&$(0.5,0.5)$		\\ \hline  
slow disk		&1023	&51.2		&{$1.8e-01$}	&{$4$}		&43	&$(0.5,0.25)$		\\ \hline
fault, left source	&1023	&51.2		&{$1.1e-01$}	&{$4$}		&48	&$(0.25,0.5)$		\\ \hline  
fault, right source	&1023	&51.2		&{$2.2e-01$}	&{$4$}		&48	&$(0.75,0.5)$		\\ \hline
diagonal fault		&1023	&51.2		&{$2.6e-01$}	&{$256$}	&101	&$(0.5,0.5)$		\\ \hline
periodic medium		&319	&6		&{$1.0e-01$}	&{$1280$}	&792	&$(0.5,0.5)$		\\ \hline
\end{tabular}
\end{center} 
\caption{For each medium considered, we show the parameters $n$ and $\omega/2\pi$, along with the resulting discretization error caused by the Finite Difference (FD error) formulation. We also show the width $w$ of the PML needed, in number of points, to obtain an error caused by the PML of less than $1e-1$. Finally, we show the total number of basis matrices ($P$) needed to probe the entire DtN map with an accuracy of about $1e-1$ as found in Section \protect\ref{sec:tests}.} 
\label{FDPMLerr} 
\end{table}

In order to get a point of reference for the accuracy benchmarks, and for small problems only, the actual matrix $D$ is computed explicitly by solving the exterior problem $4N$ times using the standard basis as Dirichlet boundary conditions (in general, one will only have access to a black-box that would output the product of $D$ with some input vector). Then, we extract a block $M$ of $D$, corresponding to the interactions between two sides of $\partial \Omega$. We note that some blocks in $D$ are the same up to transpositions or flips (inverting the order of columns or rows). We call the number of copies of a block $M$, up to transpositions or flips, its multiplicity, and write it as $m(M)$. As mentioned earlier, only the distinct blocks of $D$ need to be probed. Once we have chosen a block $M$, we may calculate the inner products $x_j = \< B_j, M\>$, which are the ideal coefficients in the expansion of $M$ (this is true since the $B_j$'s have been orthonormalized). Let $M_p =\sum_{j=1}^p x_j B_j$. We %
define the $p$-term approximation error, for the block $M$, to be the quantity
\begin{equation}\label{apperr}
\sqrt{m(M)} \frac{\|M-M_p\|_F}{\|D\|_F}
\end{equation}
using the Frobenius norm. Because the blocks on the diagonal of $D$ have a singularity, their Frobenius norm can be a few orders of magnitude greater than that of other blocks, and so it is more important to approximate those well. This is why we consider the error relative to $D$, not to the block $M$. Also, we multiply by the square root of the multiplicity of the block to give us a better idea of how big the total error $\|D-\tilde{D}\|_F/\|D\|_F$ will be. For brevity, we shall refer to (\ref{apperr}) simply as the approximation error when it is clear from the context what $M$, $D$, $p$ and the $B_i$'s are.

Then, using matrix probing, we will recover a coefficient vector $\tilde{x}$ close to $x$, which should give us an approximate $\tilde{M}=\sum_{j=1}^p \tilde{x}_j B_j$, which should be close to $M$ if our choice of basis matrices is relevant. We now define the probing error (which depends on $q$ and $z$), for the block $M$, to be the quantity
\begin{equation}\label{acterr}
\sqrt{m(M)}\frac{\|M-\tilde{M}\|_F}{\|D\|_F}. 
\end{equation}
Again, for brevity, we refer to (\ref{acterr}) as the probing error when other parameters are clear from the context.

We shall present results on the approximation and probing errors for various media, along with related condition numbers, and then we shall verify that using an approximate $\tilde{D}$ (constructed from approximate $\tilde{M}$'s for each block $M$ in $D$) does not affect the accuracy of the new solution to the Helmholtz equation.

\subsection{Probing tests}\label{sec:tests}

As we saw in Section \ref{sec:probe}, randomness plays a role in the value of $\mbox{cond}(\Psi_z)$ and of the probing error. Hence, whenever we show plots for those quantities in this section, we have done 10 trials for each value of $q$ used. The error bars show the minimum and maximum of the quantity over the 10 trials, and the line is plotted through the average value over the 10 trials. As expected, we will see in all experiments that increasing $q$ gives a better conditioning, and consequently a better accuracy and smaller failure probability. The following probing results will then be used in Section \ref{sec:insolver} to solve the Helmholtz equation.

\subsubsection{Constant medium}
\begin{figure}[ht]
\begin{minipage}[t]{0.48\linewidth}
\includegraphics[scale=.5]{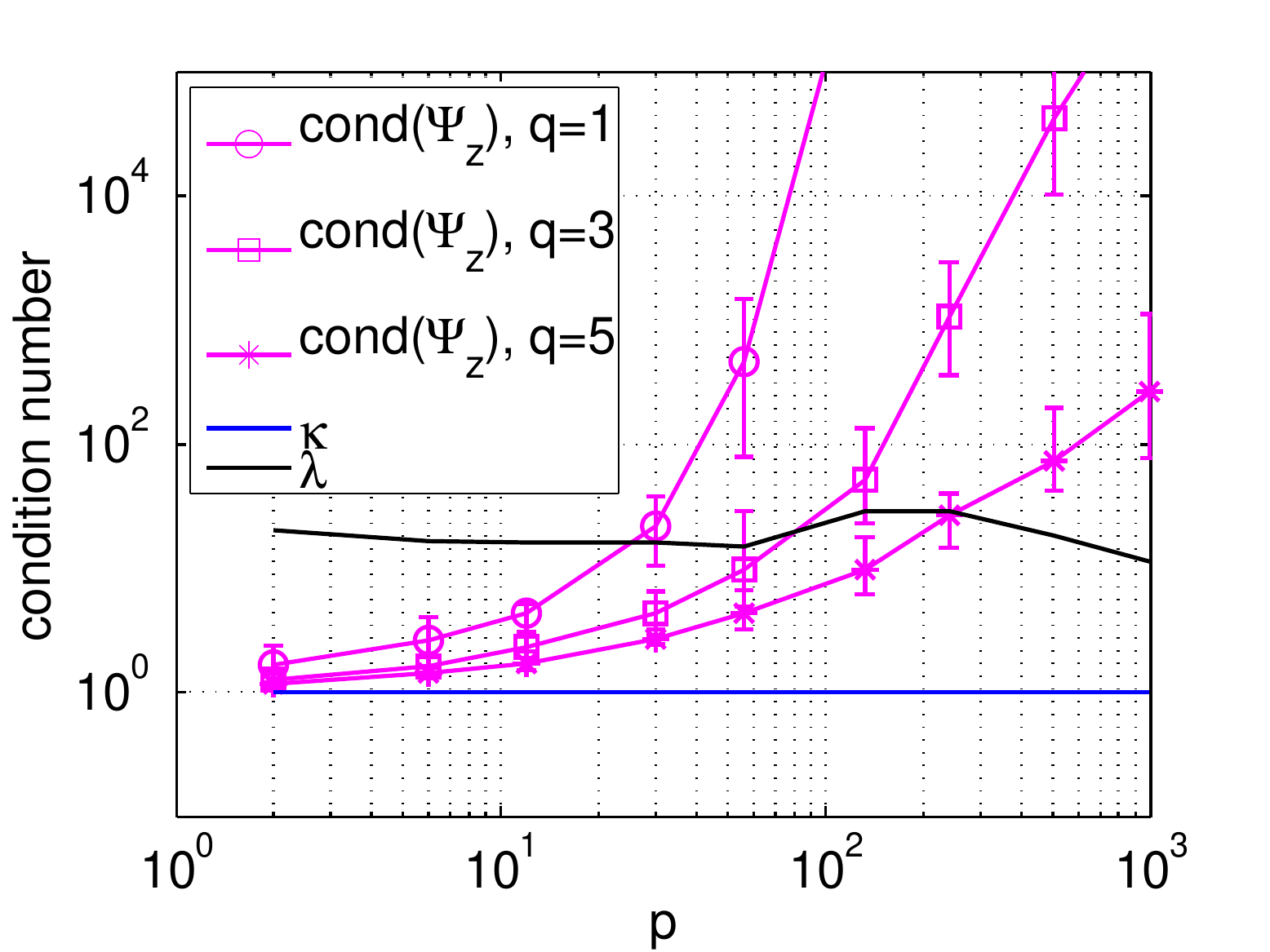}
\caption{Condition numbers for the $(1,1)$ block, $c\equiv 1$.}
\label{cond11_1024_c1}
\end{minipage}
\hspace{0.1cm}
\begin{minipage}[t]{0.48\linewidth}
\includegraphics[scale=.5]{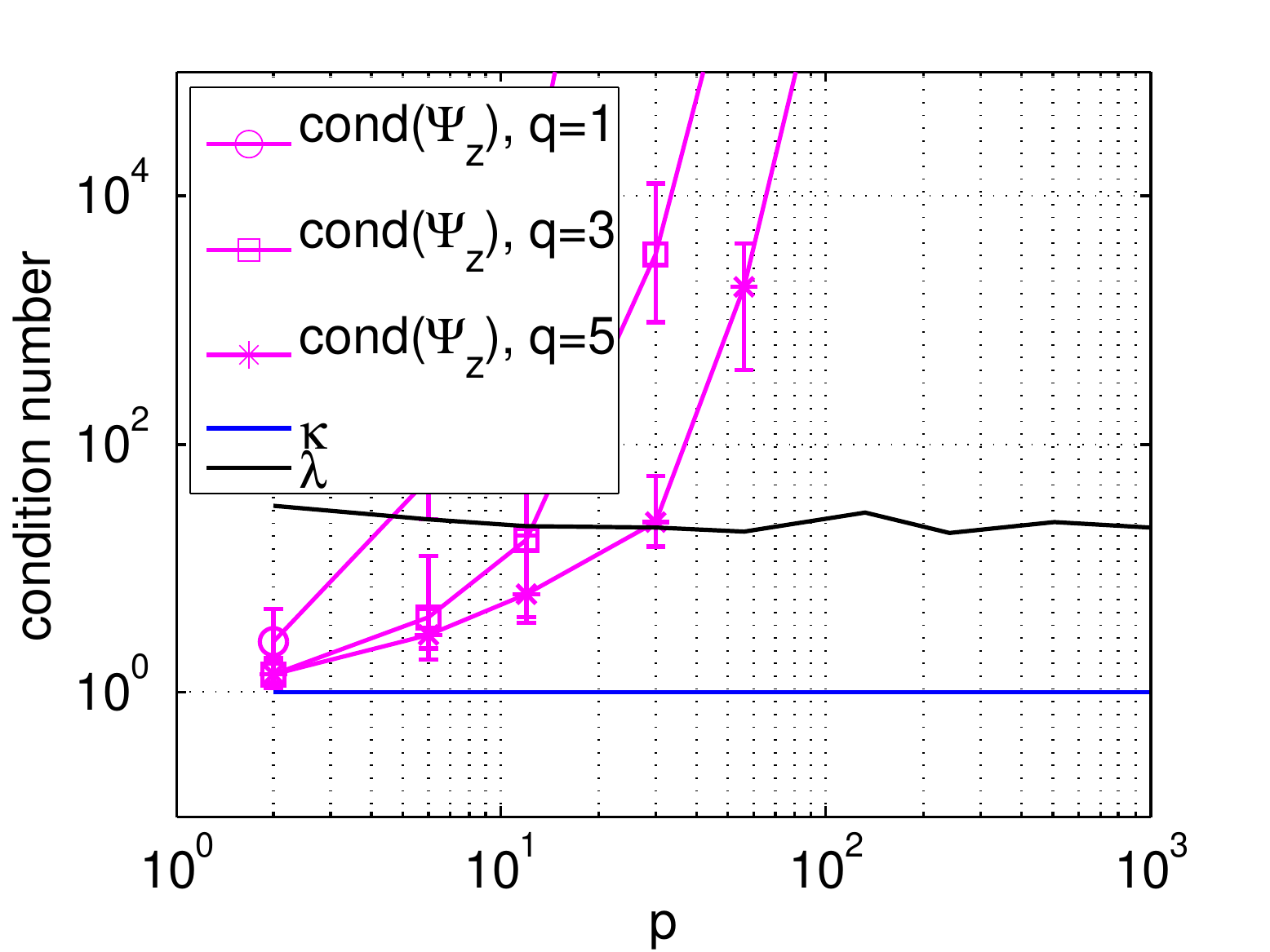}
\caption{Condition numbers for the $(2,1)$ block, $c\equiv 1$.}
\label{cond21_1024_c1}
\end{minipage}
\end{figure}

For a constant medium, $c \equiv 1$, we have three blocks with the following multiplicities: $m((1,1))=4$ (same edge), $m((2,1))=8$ (neighboring edges), and $m((3,1))=4$ (opposite edges). Note that we do not present results for the $(3,1)$ block: this block has negligible Frobenius norm\footnote{We can use probing with $q=1$ and a single basis matrix (a constant multiplied by the correct oscillations) and have a probing error of less than $10^{-6}$ for that block.} compared to $D$.  First, let us look at the conditioning for blocks $(1,1)$ and $(2,1)$. Figures \ref{cond11_1024_c1} and \ref{cond21_1024_c1} show the three relevant conditioning quantities: $\kappa$, $\lambda$ and $\mbox{cond}(\Psi_z)$ for each block. As expected, $\kappa=1$ because we orthogonalize the basis functions. Also, we see that $\lambda$ does not grow very much as $p$ increases, it remains on the order of 10. As for $\mbox{cond}(\Psi_z)$, it increases as $p$ increases for a fixed $q$ and $n$, as expected. This will affect probing in terms of the failure probability (the odds that the matrix $\Psi_z$ is far from the expected value) and accuracy (taking the pseudo-inverse will introduce larger errors in $x$). We notice these two phenomena in Figure \ref{erb1023_c1}, where we show the approximation and probing errors in probing the $(1,1)$ block for various $p$, using different $q$ and making 10 tests for each $q$ value as explained previsouly. As expected, as $p$ increases, the variations between trials get larger. Also, the probing error, always larger than the approximation error, becomes farther and farther away from the approximation error. Comparing Figure \ref{erb1023_c1} with Table \ref{c1solve} of the next section, we see that in Table \ref{c1solve} we are able to achieve higher accuracies. This is because we use the first two traveltimes (so four different types of oscillations, as explained in Section \ref{sec:basis}) to obtain those higher accuracies. But we do not use four types of oscillations for lower accuracies because this demands a larger number of basis matrices $p$ and of solves $q$ for the same error level.

\subsubsection{The waveguide}
\begin{figure}[ht]
\begin{minipage}[t]{0.48\linewidth}
\includegraphics[scale=.5]{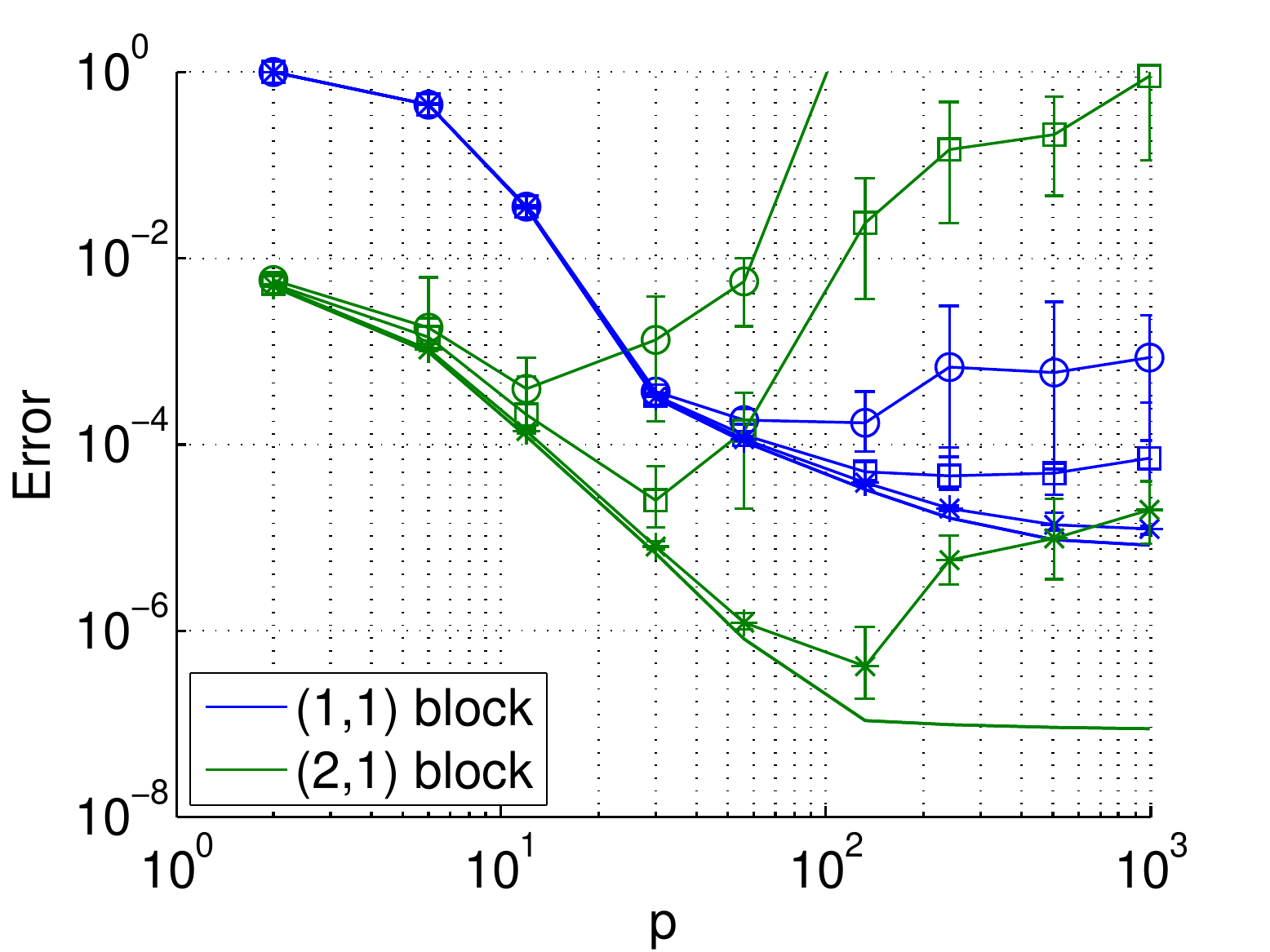}
\caption{Approximation error (line) and probing error (with markers) for the blocks of $D$, $c\equiv 1$. Circles are for $q=3$, squares for $q=5$, stars for $q=10$. }
\label{erb1023_c1}
\end{minipage}
\hspace{0.1cm}
\begin{minipage}[t]{0.48\linewidth}
\includegraphics[scale=.5]{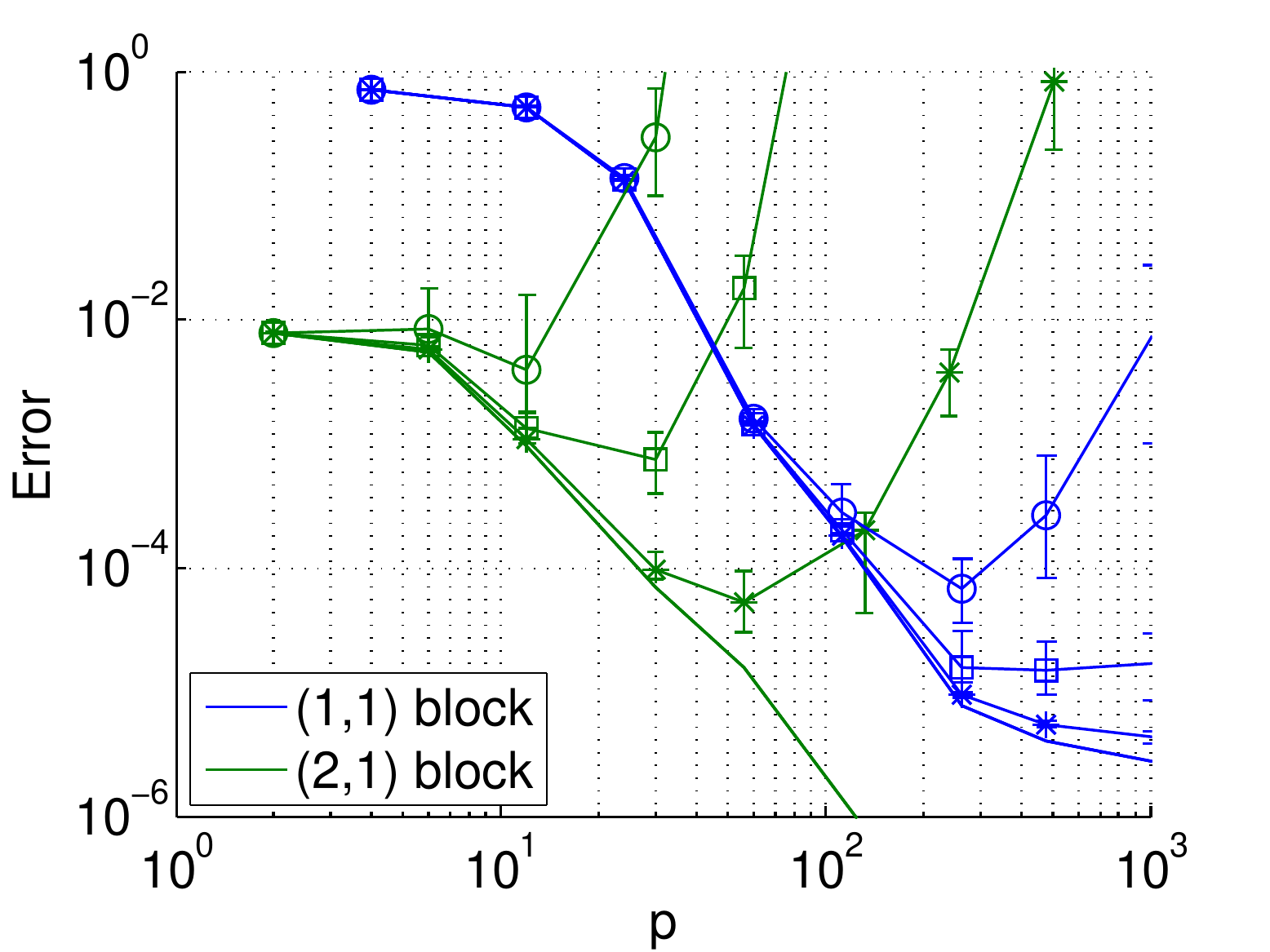}
\caption{Approximation error (line) and probing error (with markers) for the blocks of $D$, $c$ is the waveguide. Circles are for $q=3$, squares for $q=5$, stars for $q=10$.}
\label{erb1023_c2}
\end{minipage}
\end{figure}

For a waveguide as a velocity field, we have more blocks compared to the constant medium case, with different multiplicities: $m((1,1))=2$, $m((2,2))=2$, $m((2,1))=8$, $m((3,1))=2$, $m((4,2))=2$. Note that block $(2,2)$ will be easier to probe than block $(1,1)$ since the medium is smoother on that interface. Also, we can probe blocks $(3,1)$ and $(4,2)$ with $q=1$, $p=2$ and have a probing error less than $10^{-7}$. Hence we only show results for the probing and approximation errors of blocks $(1,1)$ and $(2,1)$, in Figure \ref{erb1023_c2}. Results for using probing in a solver can be found in Section \ref{sec:insolver}.

\subsubsection{The slow disk}
\begin{figure}[ht]
\begin{minipage}[t]{0.48\linewidth}
\includegraphics[scale=.5]{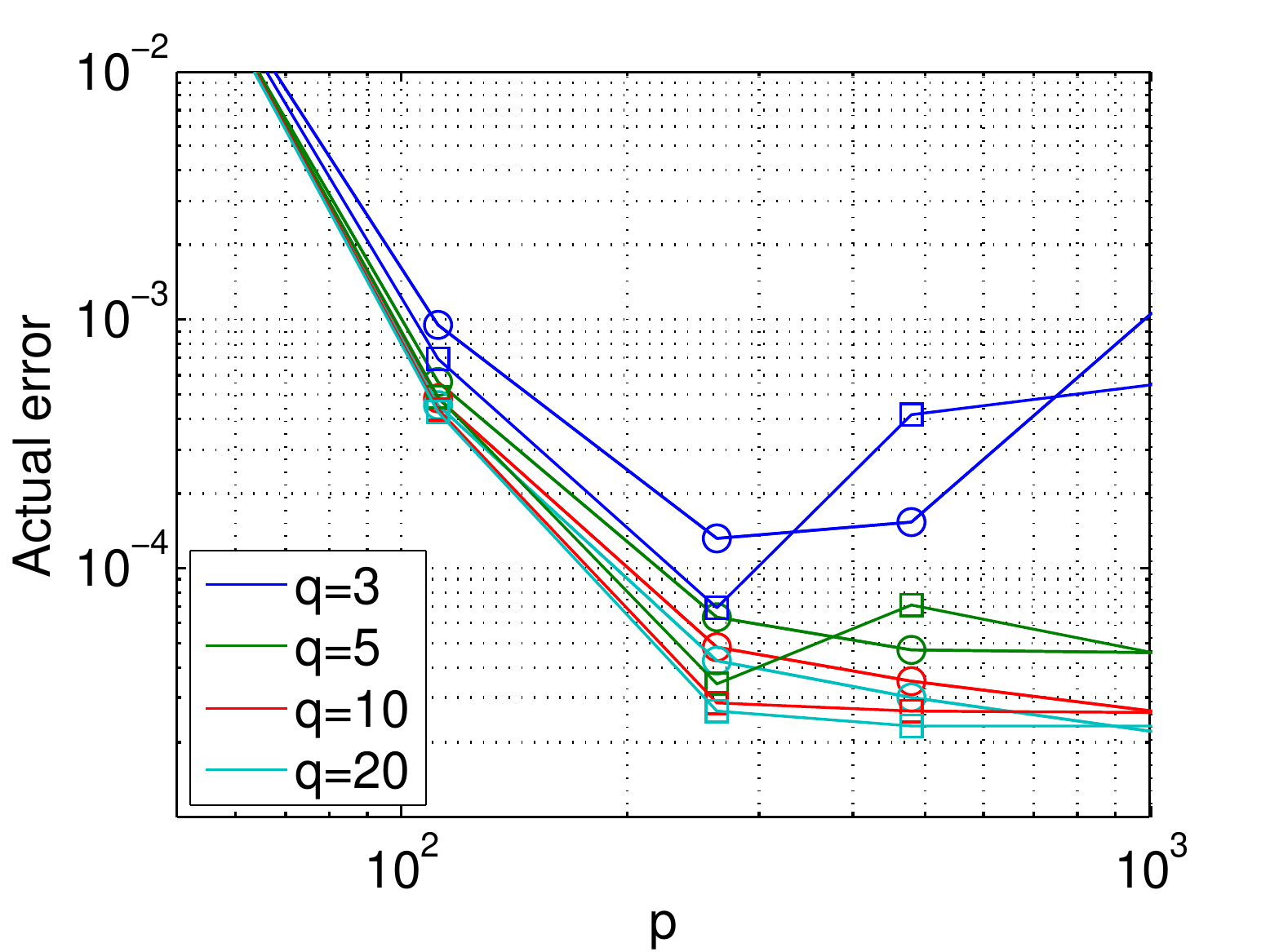}
\caption{Approximation error for the $(1,1)$ blocks of $D$, $c$ is the slowness disk, comparing the use of the normal traveltime (circles) to the fast marching traveltime (squares).}
\label{c5fastvsnorm1}
\end{minipage}
\hspace{0.1cm}
\begin{minipage}[t]{0.48\linewidth}
\includegraphics[scale=.5]{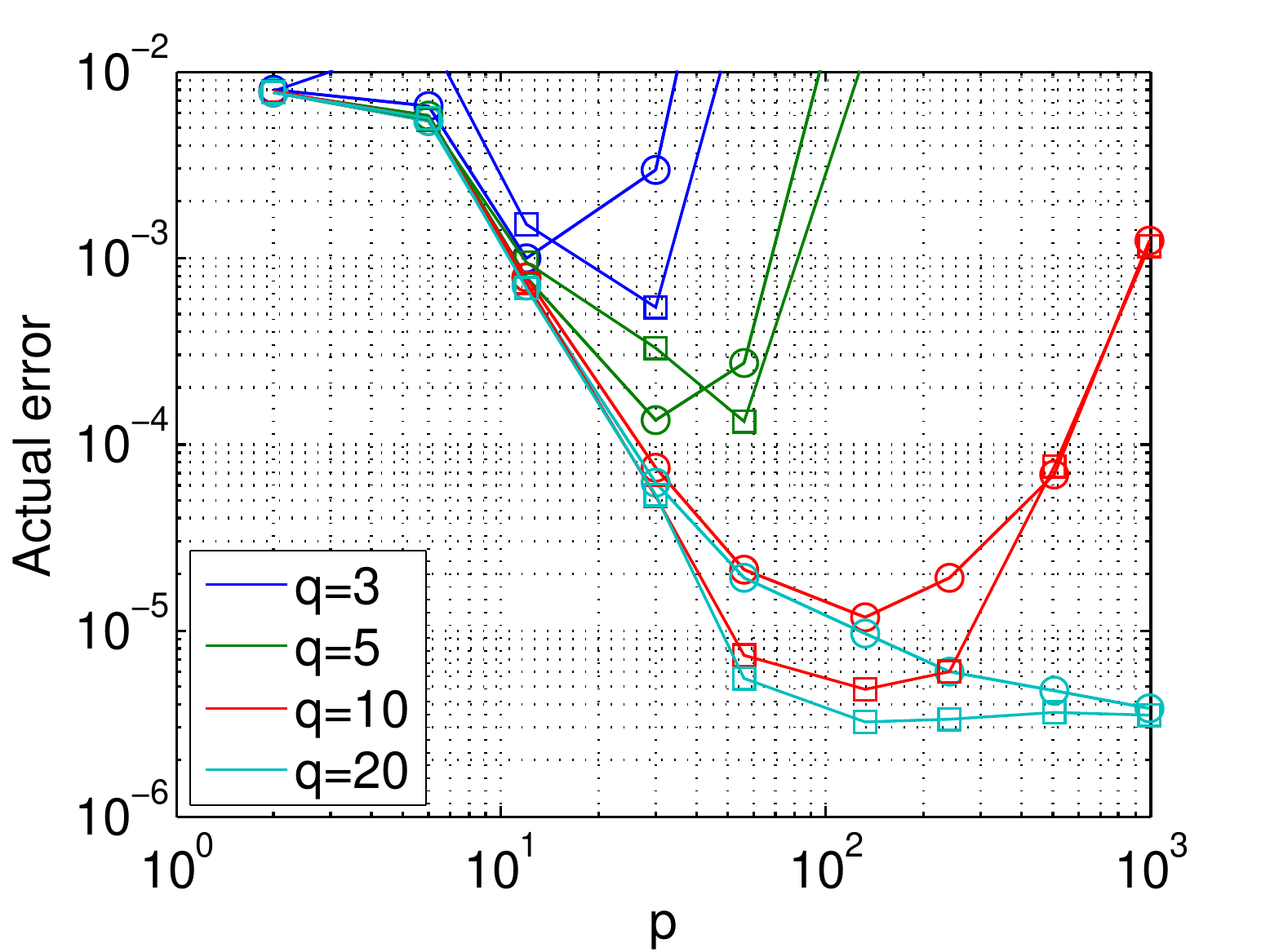}
\caption{Approximation error for the $(2,1)$ blocks of $D$, $c$ is the slowness disk, comparing the use of the normal traveltime (circles) to the fast marching traveltime (squares).}
\label{c5fastvsnorm2}
\end{minipage}
\end{figure}

Next, we consider the slow disk. Here, we have a choice to make for the traveltime upon which the oscillations depend. We may consider interfacial waves, travelling in straight line segments along $\partial \Omega$, with traveltime $\tau$. There is also the first arrival time of body waves, $\tau_f$, which for some points on $\partial \Omega$ involve taking a path that goes away from $\partial \Omega$, into the exterior where $c$ is higher, and back towards $\partial \Omega$. We have approximated this $\tau_f$ using the fast marching method of Sethian \cite{sethart}. For this example, it turns out that using either $\tau$ or $\tau_f$ to obtain oscilllations in our basis matrices does not significantly alter the probing accuracy or conditioning, although it does seem that, for higher accuracies at least, the fast marching traveltime makes convergence slightly faster. Figures \ref{c5fastvsnorm1} and \ref{c5fastvsnorm2} demonstrate this for blocks $(1,1)$ and $(2,1)$ respectively. We omit plots of the probing and approximation errors, and refer the reader to Section \ref{sec:insolver} for final probing results and using those in a solver.

\subsubsection{The vertical fault}
\begin{figure}[ht]
\begin{minipage}[t]{0.48\linewidth}
\includegraphics[scale=.5]{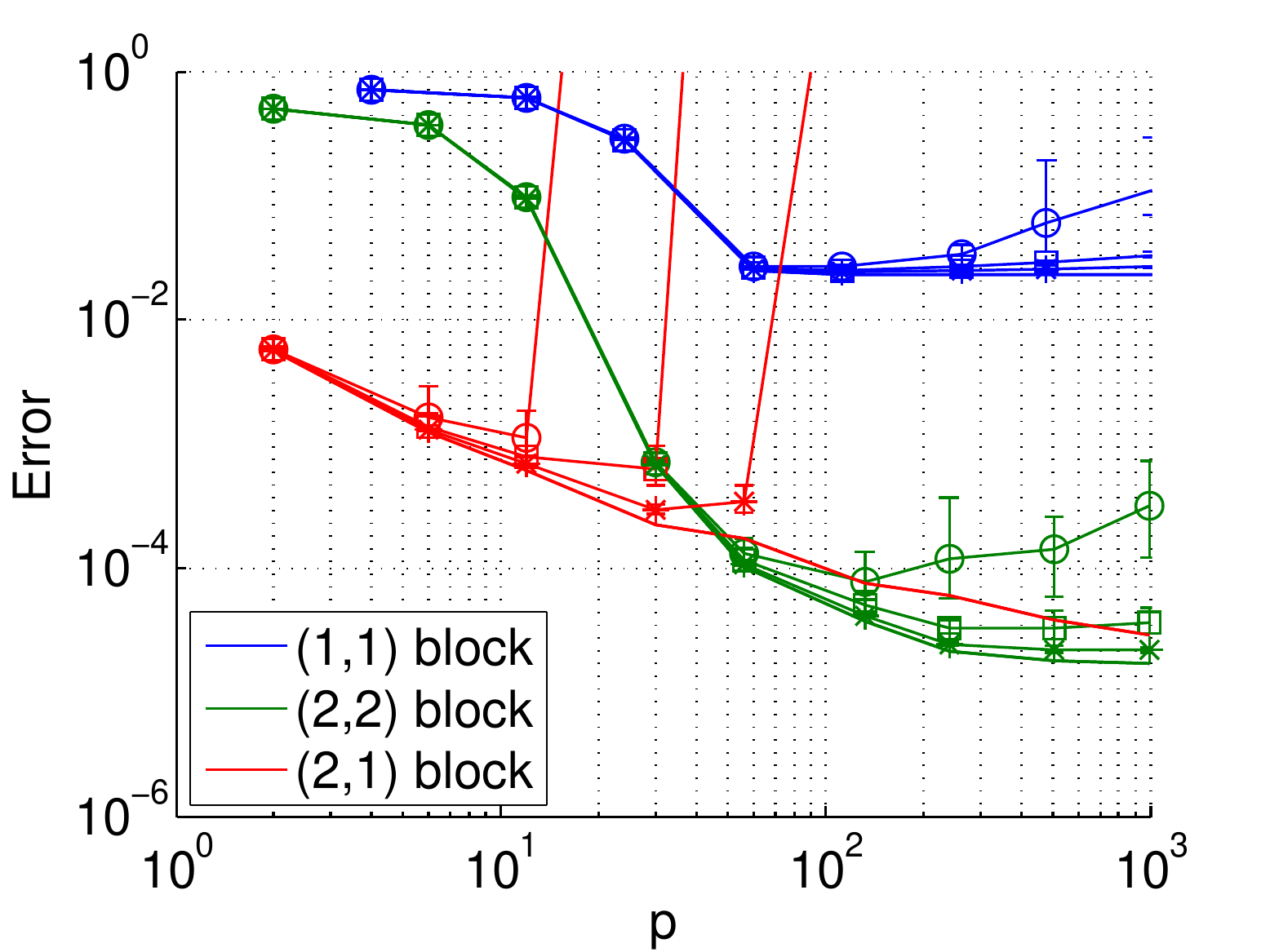}
\caption{Approximation error (line) and probing error (with markers) for the blocks of $D$, $c$ is the fault. Circles are for $q=3$, squares for $q=5$, stars for $q=10$.}
\label{erb1023_c16}
\end{minipage}
\hspace{0.1cm}
\begin{minipage}[t]{0.48\linewidth}
\includegraphics[scale=.5]{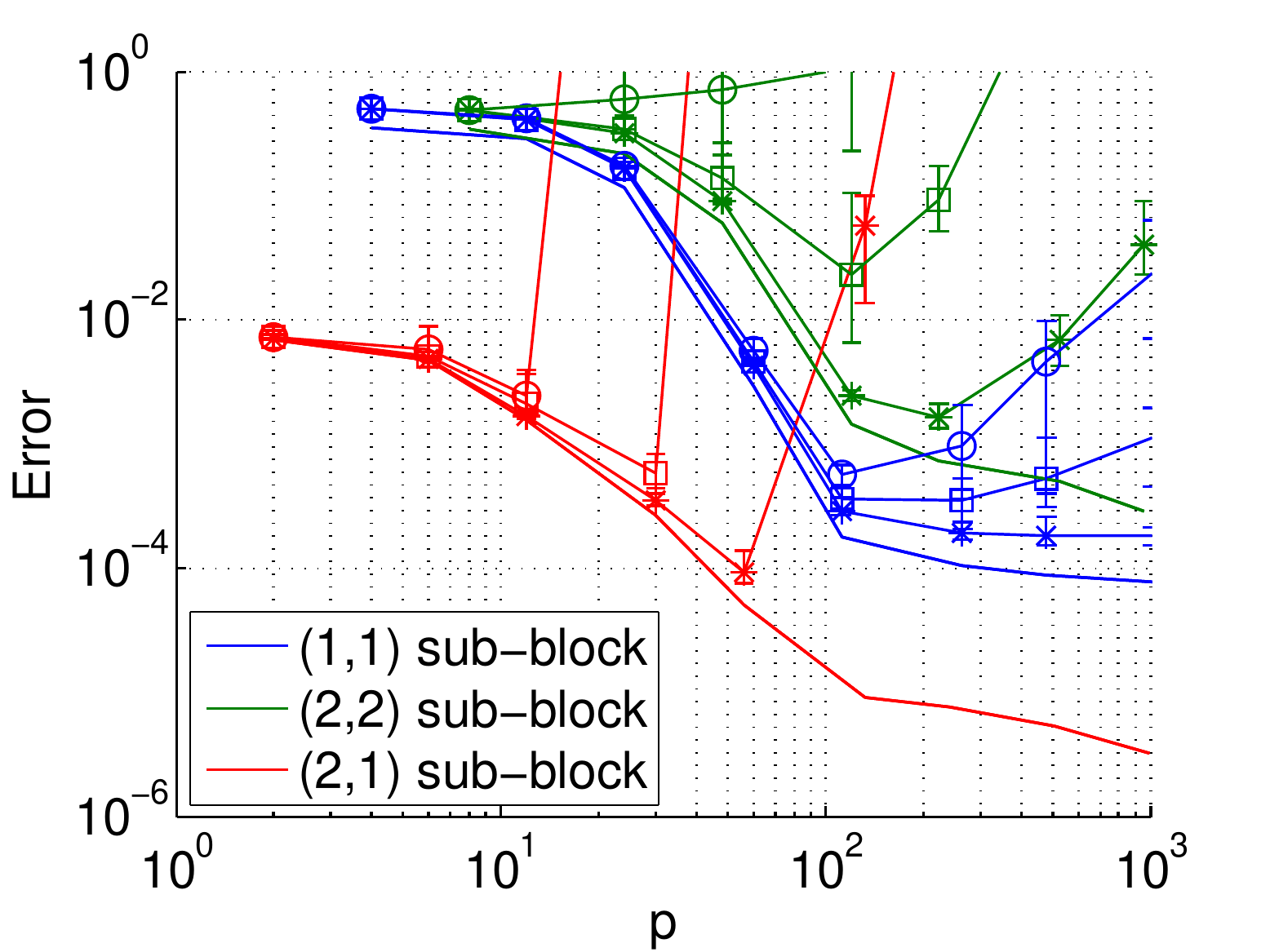}
\caption{Approximation error (line) and probing error (with markers) for the sub-blocks of the $(1,1)$ block of $D$, $c$ is the fault. Circles are for $q=3$, squares for $q=5$, stars for $q=10$.}
\label{erb1023_c16_sub}
\end{minipage}
\end{figure}

Next, we look at the case of the medium $c$ which has a vertical fault. We note that this case is harder because some of the blocks will have themselves a 2 by 2 or 1 by 2 structure caused by the discontinuity in the medium. Ideally, as we shall see, each sub-block should be probed separately. There are 7 distinct blocks, with different multiplicities: $m((1,1))=2$, $m((2,2))=1$, $m((4,4))=1$, $m((2,1))=4$, $m((4,1))=4$, $m((3,1))=2$, $m((4,2))=2$. Blocks $(2,2)$ and $(4,4)$ are easier to probe than block $(1,1)$ because they do not exhibit a sub-structure. Also, since the velocity is smaller on the right side of the fault, the frequency there is higher, which means that blocks involving side 2 are slightly harder to probe than those involving side 4. Hence we first present results for the blocks $(1,1)$, $(2,2)$ and $(2,1)$ of $D$. In Figure \ref{erb1023_c16} we see the approximation and probing errors for those blocks. Then, in Figure \ref{erb1023_c16_sub}, we present results for the errors related to probing the 3 distinct sub-blocks of the $(1,1)$ block of $D$. We can see that probing the $(1,1)$ block by sub-blocks helps achieve greater accuracy. We could have split other blocks too to improve the accuracy of their probing (for example, block $(2,1)$ has a 1 by 2 structure because side 1 has a discontinuity in $c$) but the accuracy of the overall DtN map was still limited by the accuracy of probing the $(1,1)$ block, so we do not show results for other splittings.

\subsubsection{The diagonal fault}
Now, we look at the case of the medium $c$ which has a diagonal fault. Again, some of the blocks will have themselves a 2 by 2 or 1 by 2 structure. There are 6 distinct blocks, with different multiplicities: $m((1,1))=2$, $m((2,2))=2$, $m((2,1))=4$, $m((4,1))=2$, $m((3,2))=2$, $m((3,1))=4$. Again, we split up block $(1,1)$ in 4 sub-blocks and probe each of those sub-blocks separately for greater accuracy, but do not split other blocks. We then use two traveltimes for the $(2,2)$ sub-block of block $(1,1)$. Using as the second arrival time the geometrical traveltime consisting of leaving the boundary and bouncing off the fault, as mentioned in Section \ref{sec:tt}, allowed us to increase accuracy by an order of magnitude compared to using only the first arrival traveltime, or compared to using as a second arrival time the usual bounce off the corner (or here, bounce off the fault where it meets $\delta \Omega$). We omit plots of the probing and approximation errors, and refer the reader to Section \ref{sec:insolver} for final probing results and using those in a solver.

\subsubsection{The periodic medium}

Finally, we look at the case of the periodic medium presented earlier. There are 3 distinct blocks, with different multiplicities: $m((1,1))=4$, $m((2,1))=8$, $m((3,1))=4$. We expect the corresponding DtN map to be harder to probe because its structure will reflect that of the medium, i.e. it will exhibit sharp transitions at points corresponding to sharp transitions in $c$ (similarly as with the faults). First, we notice that, in all the previous mediums we tried, plotting the norm of the anti-diagonal entries of diagonal blocks (or sub-blocks for the faults) shows a rather smooth decay away from the diagonal. However, that is not the case for the periodic medium: it looks like there is decay away from the diagonal, but variations from that decay can be of relative order 1. This prevents our usual strategy, using basis matrices containing terms that decay away from the diagonal such as $(h+|x-y|)^{-j_1/\alpha}$, from working adequately. Instead, we use polynomials along anti-diagonals, as well as polynomials along diagonals as we previously did.

It is known that solutions to the Helmholtz equation in a periodic medium are Bloch waves with a particular structure \cite{JohnsonPhot}. However, using that structure in the basis matrices is not robust. Indeed, using a Bloch wave structure did not succeed very well, probably because our discretization was not accurate enough and so $D$ exhibited that structure only to a very rough degree. Hence we did not use Bloch waves for probing the periodic medium.

For this reason, we tried basis matrices with no oscillations, but with polynomials in both directions as explained previously, and obtained the results of Section \ref{sec:insolver}. Now that we have probed the DtN map and obtained compressed blocks to form an approximation $\tilde{D}$ of $D$, we may use this $\tilde{D}$ in a Helmholtz solver as an absorbing boundary condition.

\subsection{Using the probed $D$ into a Helmholtz solver}\label{sec:insolver}

\begin{figure}[ht]
\begin{minipage}[t]{0.45\linewidth}
\includegraphics[scale=.45]{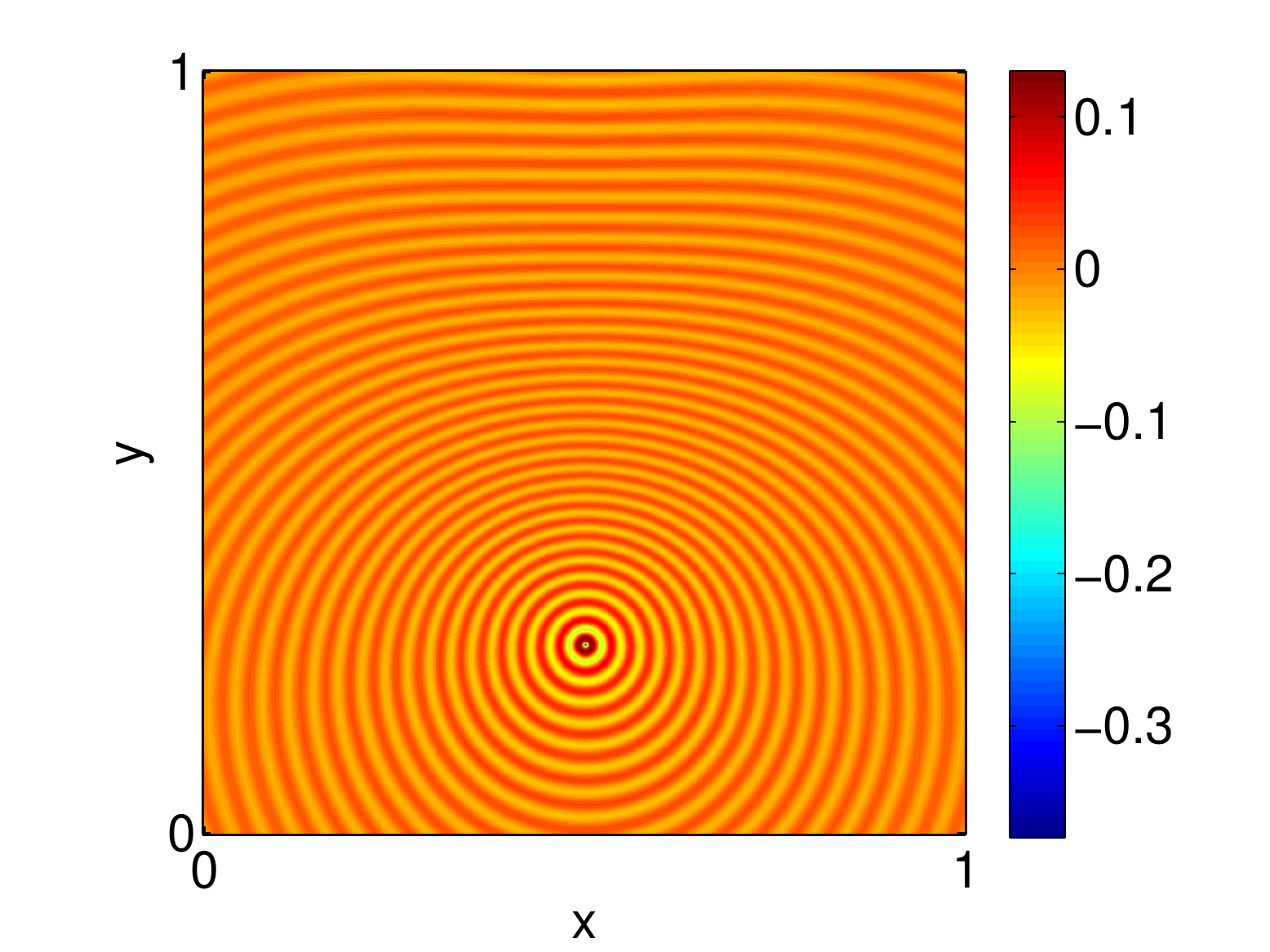}
\caption{Real part of the solution, $c$ is the slow disk.}
\label{solc5}
\end{minipage}
\hspace{1cm}
\begin{minipage}[t]{0.45\linewidth}
\includegraphics[scale=.45]{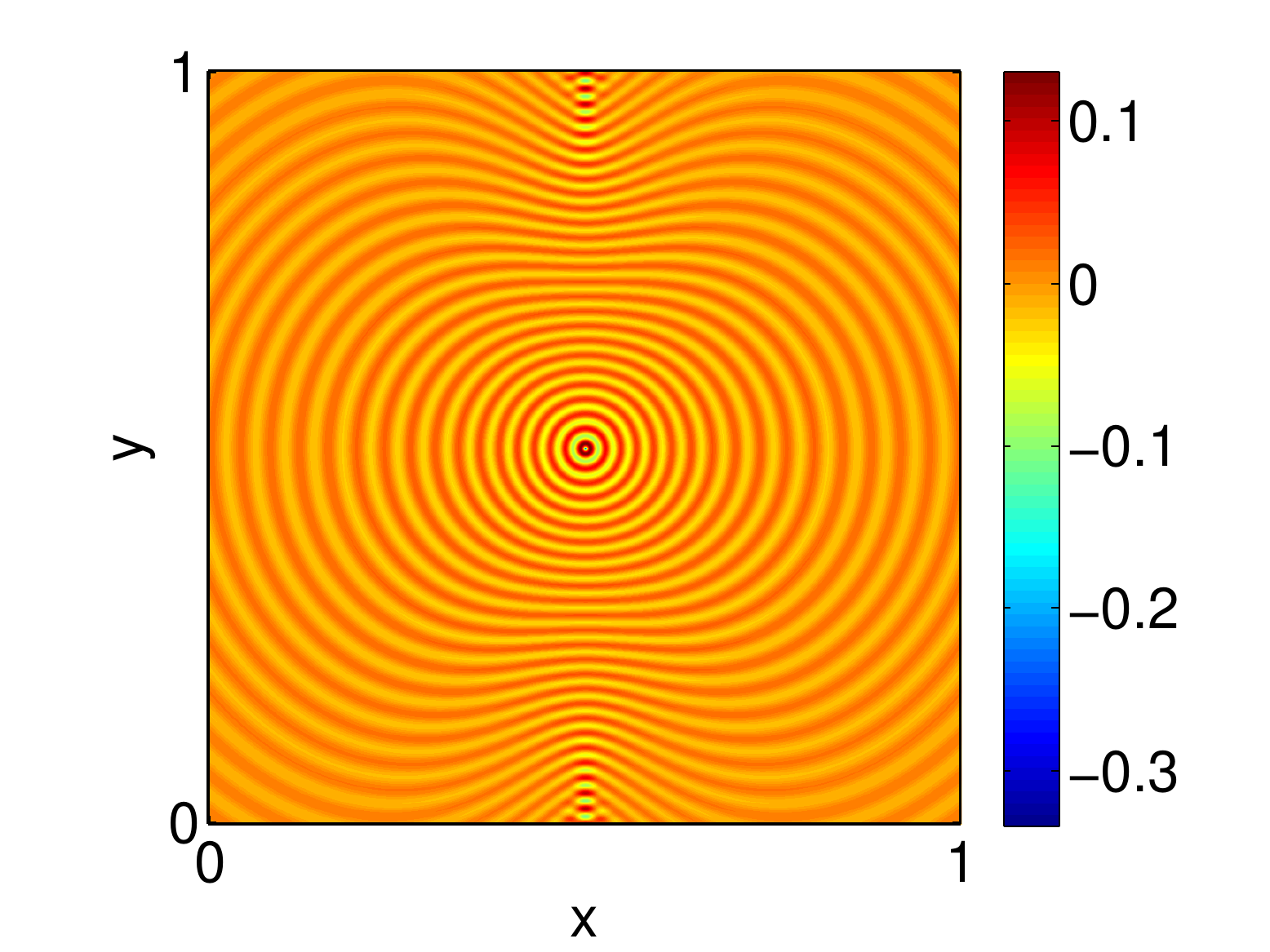}
\caption{Real part of the solution, $c$ is the waveguide.}
\label{solc3}
\end{minipage}
\end{figure}

\begin{figure}[ht]
\begin{minipage}[t]{0.45\linewidth}
\includegraphics[scale=.45]{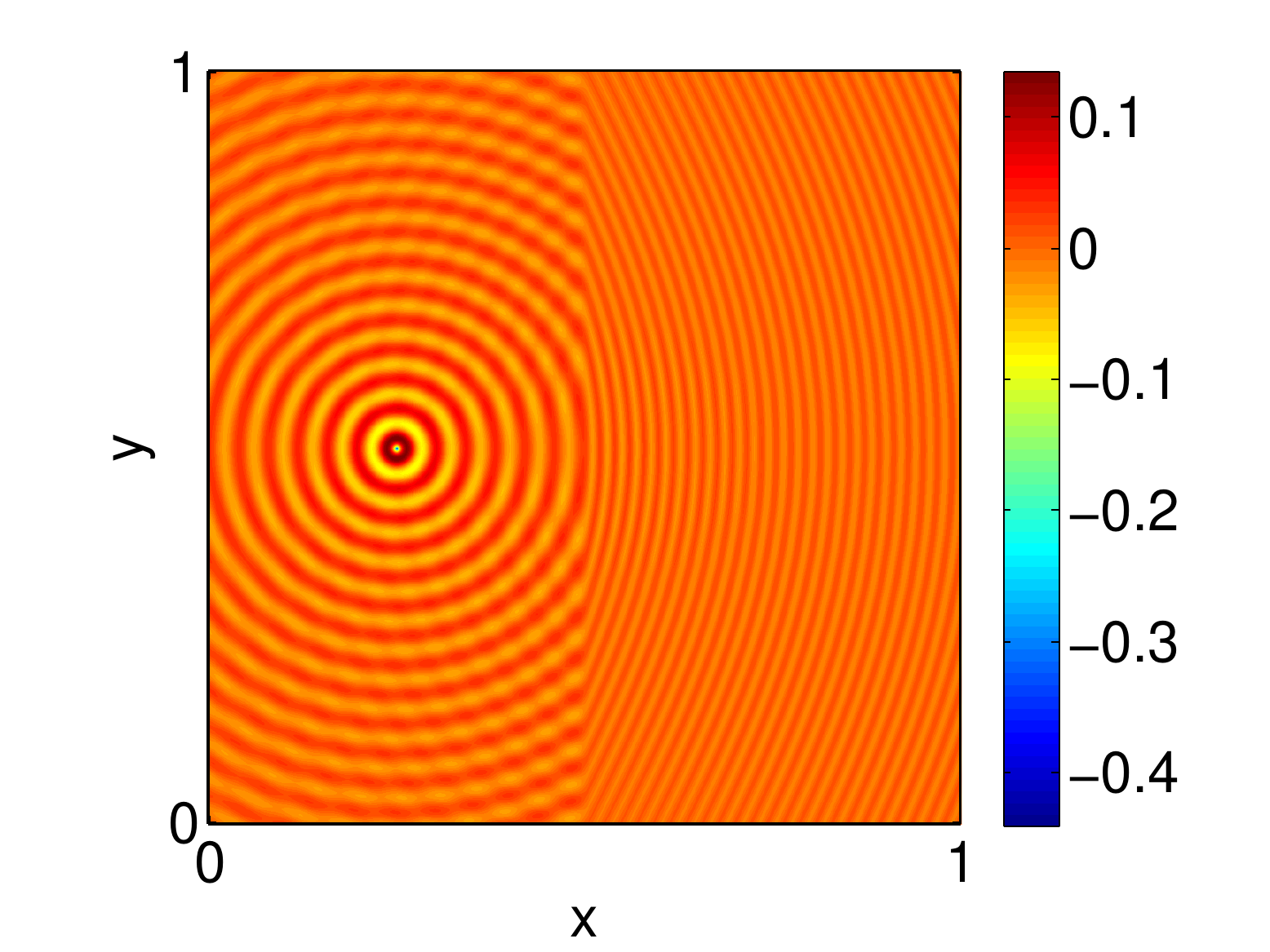}
\caption{Real part of the solution, $c$ is the vertical fault with source on the left.}
\label{solc16l}
\end{minipage}
\hspace{1cm}
\begin{minipage}[t]{0.45\linewidth}
\includegraphics[scale=.45]{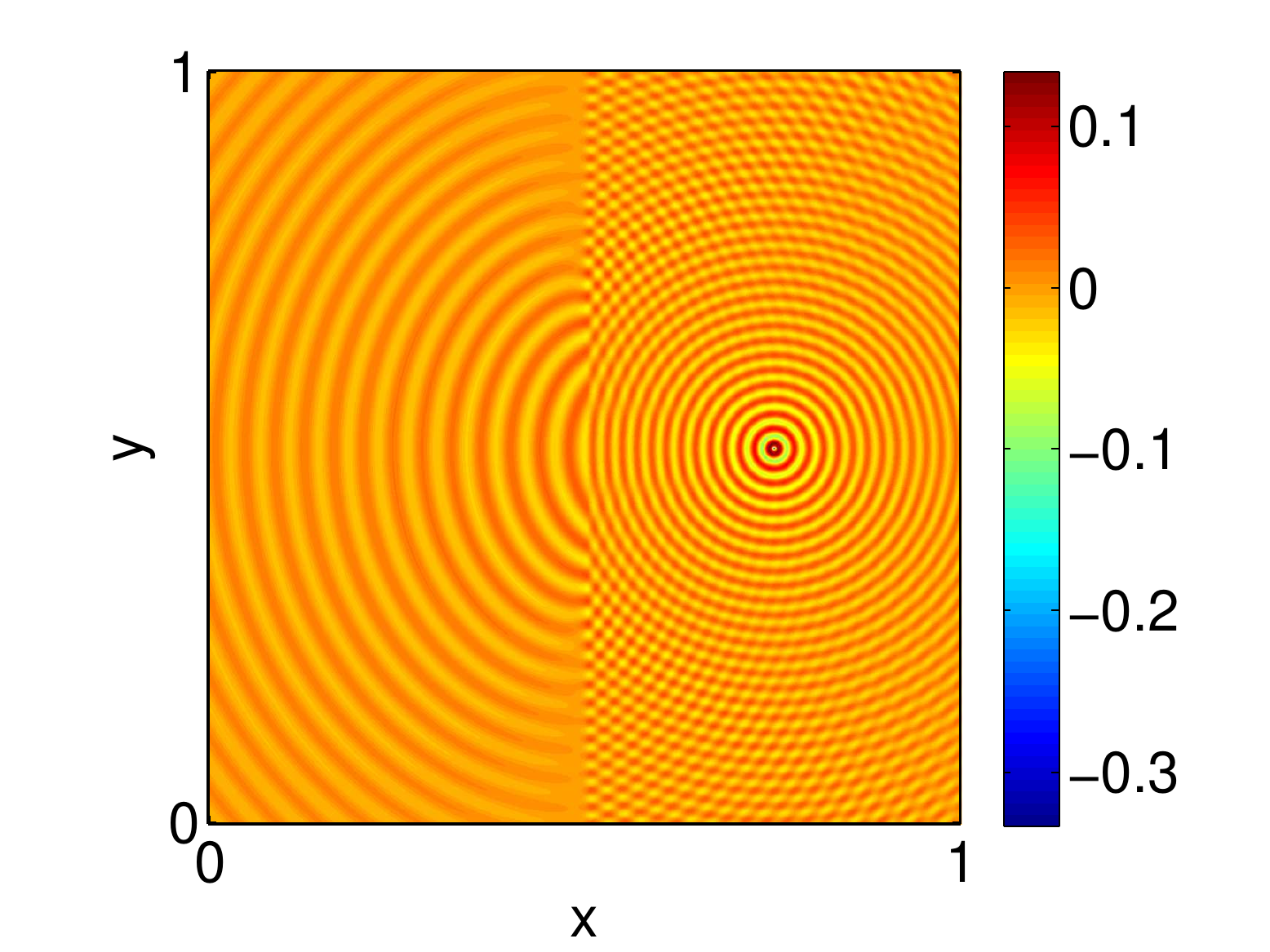}
\caption{Real part of the solution, $c$ is the vertical fault with source on the right.}
\label{solc16r}
\end{minipage}
\end{figure}

\begin{figure}[ht]
\begin{minipage}[t]{0.45\linewidth}
\includegraphics[scale=.45]{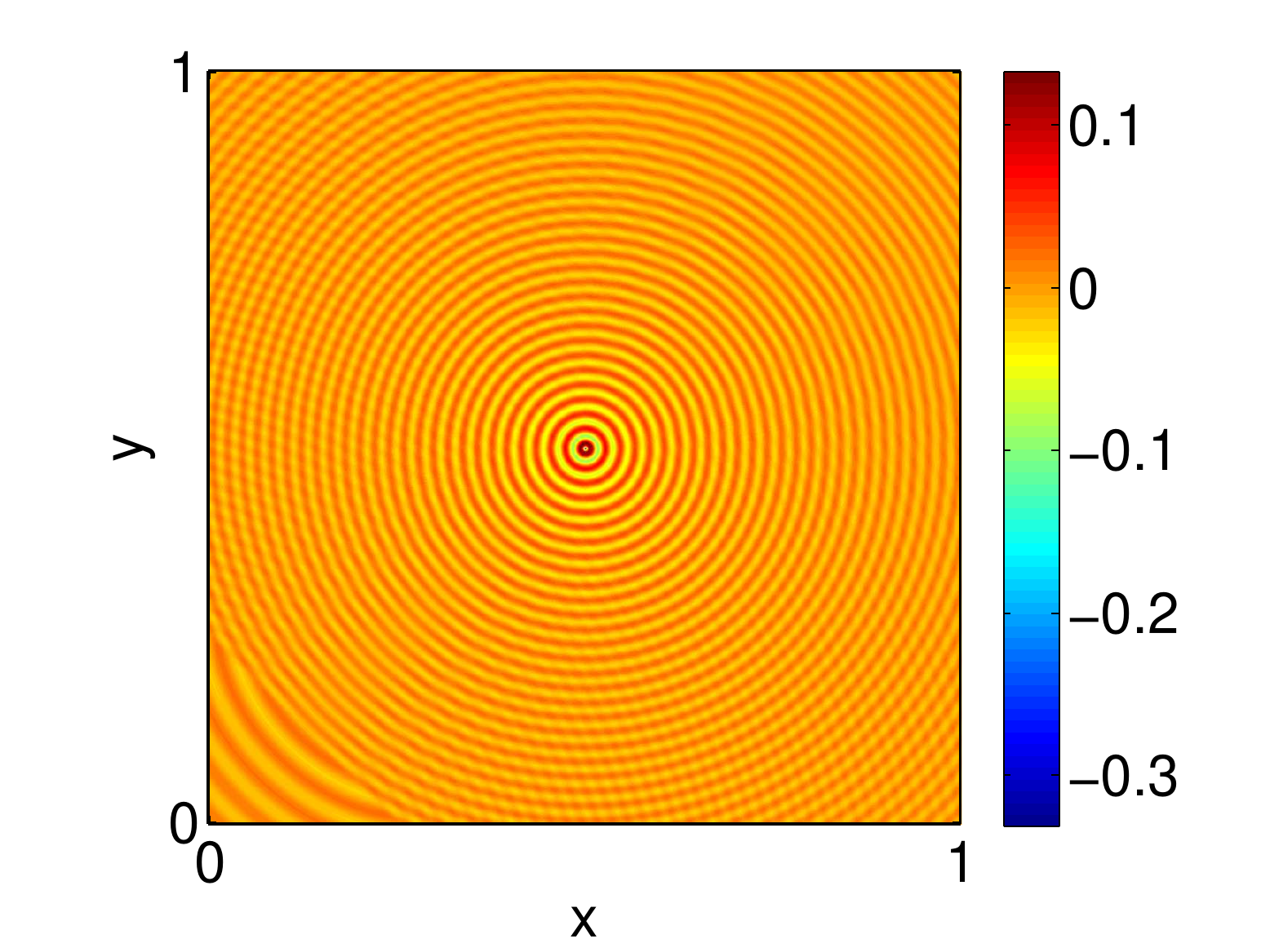}
\caption{Real part of the solution, $c$ is the diagonal fault.}
\label{solc18}
\end{minipage}
\hspace{1cm}
\begin{minipage}[t]{0.45\linewidth}
\includegraphics[scale=.45]{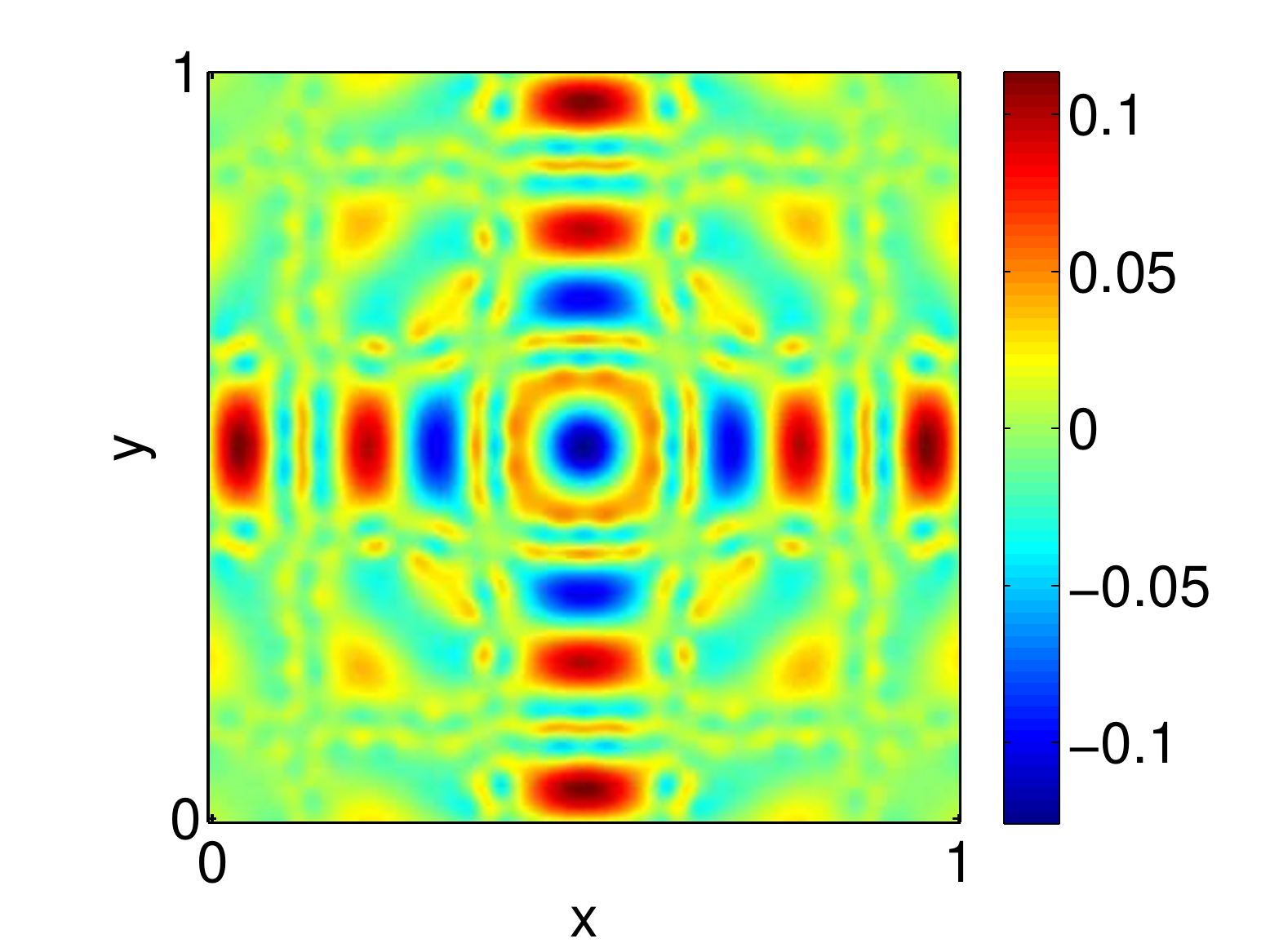}
\caption{Imaginary part of the solution, $c$ is the periodic medium.}
\label{solc33}
\end{minipage}
\end{figure}

In Figures \ref{solc5}, \ref{solc3}, \ref{solc16l}, \ref{solc16r}, \ref{solc18} and \ref{solc33} we can see the standard solutions to the Helmholtz equation on $[0,1]^2$ using a PML for the various media we consider, except for the constant medium, where the solution is well-known. We use those as our reference solutions.

We have tested the solver with the probed $\tilde{D}$ as an absorbing boundary condition with success. See Tables \ref{c1solve}, \ref{c5solve}, \ref{c3solve} and \ref{c16solve} for results corresponding to each medium. We show the number $p$ of basis matrices required for some blocks for that tolerance, the number of solves $q$ of the exterior problem for those blocks, the total number of solves $Q$, the error in $D$ and the relative error in Frobenius norm between the solution $\tilde{u}$ using $\tilde{D}$ and the solution $u$ using $D$.
As we can see from the tables, the error $\frac{\|u-\tilde{u}\|_F}{\|u\|_F}$ in the solution $u$ is no more than an order of magnitude greater than the error $\frac{\|D-\tilde{D}\|_F}{\|D\|_F}$ in the DtN map $D$.Grazing waves, which can arise when the source is close to the boundary of the computational domain, will be discussed in Section \ref{sec:graz}.
We note again that, for the constant medium, using the second arrival traveltime as well as the first for the $(1,1)$ block allowed us to achieve accuracies of 5 and 6 digits in the DtN map, which was not possible otherwise. Using a second arrival time for the cases of the faults was also useful. Those results show that probing works best when the medium $c$ is rather smooth. For non-smooth media such as a fault, it becomes harder to probe the DtN map to a good accuracy, so that the solution to the Helmholtz equation also contains more error. 

\begin{table}
\caption{$c\equiv 1$} 
\begin{center} \footnotesize
\begin{tabular}{|l|l|l|l|l|l|} \hline  
$p$ for $(1,1)$		& $p$ for $(2,1)$	& $q=Q$		& $\frac{\|D-\tilde{D}\|_F}{\|D\|_F}$ 	& $\frac{\|u-\tilde{u}\|_F}{\|u\|_F}$ \\ \hline  
{$6$}			& {$1$}			& {$1$}		& {$2.0130e-01$}			& {$3.3191e-01$} \\ \hline  
{$12$}			& {$2$}			& {$1$}		& {$9.9407e-03$}			& {$1.9767e-02$} \\ \hline  
{$20$}			& {$12$}		& {$3$}		& {$6.6869e-04$}			& {$1.5236e-03$} \\ \hline  
{$72$}			& {$20$}		& {$5$}		& {$1.0460e-04$}			& {$5.3040e-04$} \\ \hline  
{$224$}			& {$30$}		& {$10$}	& {$8.2892e-06$}			& {$9.6205e-06$} \\ \hline  
{$360$}			& {$90$}		& {$10$}	& {$7.1586e-07$}			& {$1.3044e-06$} \\ \hline  
\end{tabular}
\end{center} 
\label{c1solve} 
\end{table}

\begin{table}
\caption{$c$ is the waveguide} 
\begin{center} \footnotesize
\begin{tabular}{|l|l|l|l|l|l|l|l|l|} \hline  
$p$ for $(1,1)$	&$p$ for $(2,1)$&$q$	&$p$ for $(2,2)$&$q$  &$Q$	&$\frac{\|D-\tilde{D}\|_F}{\|D\|_F}$	&$\frac{\|u-\tilde{u}\|_F}{\|u\|_F}$ \\ \hline  
$40$		&$2$		&$1$	&$12$		&$1$  &$2$	&$9.1087e-02$				&$1.2215e-01$ \\ \hline  
$40$		&$2$		&$3$	&$20$		&$1$  &$4$	&$1.8685e-02$				&$7.6840e-02$       \\ \hline  
$60$		&$20$		&$5$	&$20$		&$3$  &$8$	&$2.0404e-03$				&$1.3322e-02$       \\ \hline  
$112$		&$30$		&$10$	&$30$		&$3$  &$13$	&$2.3622e-04$				&$1.3980e-03$      \\ \hline  
$264$		&$72$		&$20$	&$168$		&$10$ &$30$	&$1.6156e-05$				&$8.9911e-05$      \\ \hline  
$1012$		&$240$		&$20$	&$360$		&$10$ &$30$	&$3.3473e-06$				&$1.7897e-05$      \\ \hline  
\end{tabular}
\end{center} 
\label{c3solve} 
\end{table}

\begin{table}
\caption{$c$ is the slow disk} 
\begin{center} \footnotesize
\begin{tabular}{|l|l|l|l|l|l|} \hline  
$p$ for $(1,1)$	& $p$ for $(2,1)$	&$q=Q$	&$\frac{\|D-\tilde{D}\|_F}{\|D\|_F}$ 	& $\frac{\|u-\tilde{u}\|_F}{\|u\|_F}$ \\ \hline  
{$40$}		& {$2$}			&{$3$}	& {$1.0730e-01$}			& {$5.9283e-01$} \\ \hline  
{$84$}		& {$2$}			&{$3$}	& {$8.0607e-03$}			& {$4.5735e-02$} \\ \hline  
{$180$}		& {$12$}		&{$3$}	& {$1.2215e-03$}			& {$1.3204e-02$} \\ \hline  
{$264$}		& {$30$}		&{$5$}	& {$1.5073e-04$}			& {$7.5582e-04$} \\ \hline  
{$1012$}	& {$132$}		&{$20$}	& {$2.3635e-05$}			& {$1.5490e-04$} \\ \hline  
\end{tabular}
\end{center} 
\label{c5solve} 
\end{table}

\begin{table}
\caption{$c$ is the fault} 
\begin{center} \footnotesize
\begin{tabular}{|l|l|l|l|l|l|l|l|l|} \hline  
Q 	&$\frac{\|D-\tilde{D}\|_F}{\|D\|_F}$	&$\frac{\|u-\tilde{u}\|_F}{\|u\|_F}$, left source 	&$\frac{\|u-\tilde{u}\|_F}{\|u\|_F}$, right source\\ \hline  
{$5$}	&{$2.8376e-01$}				&{$6.6053e-01$} 					&{$5.5522e-01$} \\ \hline  
{$5$}	&{$8.2377e-03$}				&{$3.8294e-02$} 					&{$2.4558e-02$} \\ \hline  
{$30$}	&{$1.1793e-03$}				&{$4.0372e-03$} 					&{$2.9632e-03$} \\ \hline  
\end{tabular}
\end{center} 
\label{c16solve} 
\end{table}

\begin{table}
\caption{$c$ is the diagonal fault} 
\begin{center} \footnotesize
\begin{tabular}{|l|l|l|l|l|l|l|l|l|} \hline  
Q 	&$\frac{\|D-\tilde{D}\|_F}{\|D\|_F}$	&$\frac{\|u-\tilde{u}\|_F}{\|u\|_F}$ \\ \hline  
{$4$}	&{$1.6030e-01$}				&{$4.3117e-01$}   \\ \hline  
{$6$}	&{$1.7845e-02$}				&{$7.1500e-02$}   \\ \hline  
{$23$}	&{$4.2766e-03$}				&{$1.2429e-02$}   \\ \hline  
\end{tabular}
\end{center} 
\label{c18solve} 
\end{table}

\begin{table}
\caption{$c$ is the periodic medium} 
\begin{center} \footnotesize
\begin{tabular}{|l|l|l|l|l|l|l|l|l|} \hline  
Q 	&$\frac{\|D-\tilde{D}\|_F}{\|D\|_F}$	&$\frac{\|u-\tilde{u}\|_F}{\|u\|_F}$ \\ \hline  
{$50$}	&{$1.8087e-01$}				&{$1.7337e-01$}   \\ \hline  
{$50$}	&{$3.5714e-02$}				&{$7.1720e-02$}   \\ \hline  
{$50$}	&{$9.0505e-03$}				&{$2.0105e-02$}   \\ \hline  
\end{tabular}
\end{center} 
\label{c33solve}
\end{table}

\subsection{Grazing waves}\label{sec:graz}

\begin{figure}[ht]
\begin{center}
\includegraphics[scale=.5]{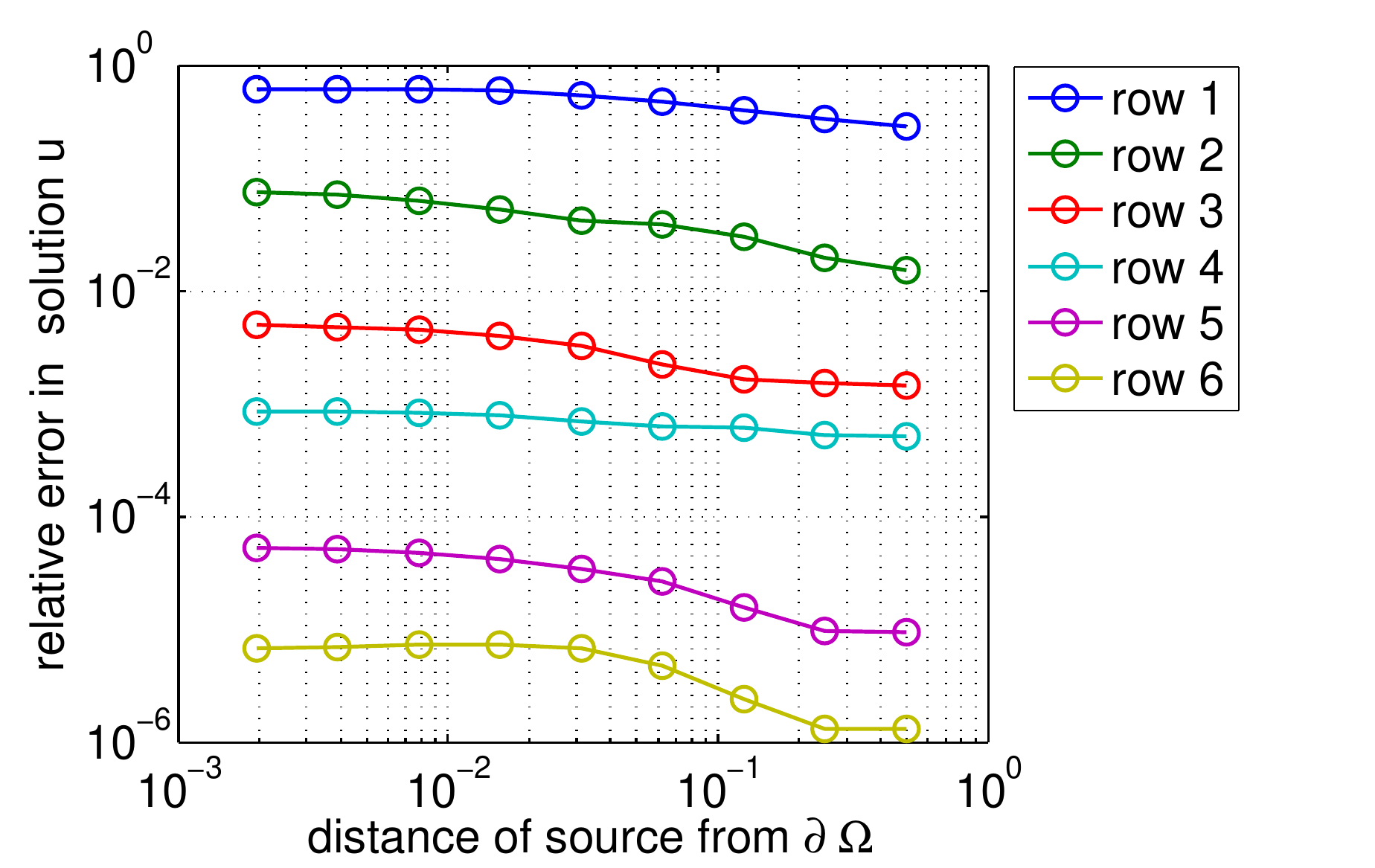}
\caption{Error in solution $u$, $c\equiv 1$, moving point source.}
\label{c1graz}
\end{center}
\end{figure}

It is well-known that ABCs often have difficulties when a source is close to a boundary of the domain, or in general when waves incident to the boundary are almost parallel to it. We wish to verify that the solution $\tilde{u}$ using the result $\tilde{D}$ of probing $D$ does not degrade as the source becomes closer and closer to some side of $\partial \Omega$. For this, we use a right hand side $f$ to the Helmholtz equation which is a point source, located at the point $(x_0,y_0)$, where $x_0=0.5$ is fixed and $y_0>0$ becomes smaller and smaller, until it is a distance $2h$ away from the boundary (the point source's stencil has width $h$, so a source at a distance $h$ from the boundary does not make sense). We see in Figure \ref{c1graz} that, for $c \equiv 1$, the solution remains quite good until the source is a distance $2h$ away from the boundary. In particular, the better the solution is for a source in the middle of the domain, the less it degrades as the source gets closer to the boundary. We obtain very similar results for the waveguide, slow disk and faults (for the vertical fault we locate the source at $(x_0,y_0)$, where $y_0=0.5$ is fixed and $x_0$ goes to $0$ or $1$). This shows that the probing process itself does not significantly affect how well grazing waves are absorbed.

\subsection{Variations of $p$ with $N$}

\begin{figure}[ht]
\begin{minipage}[t]{0.45\linewidth}
\includegraphics[scale=.45]{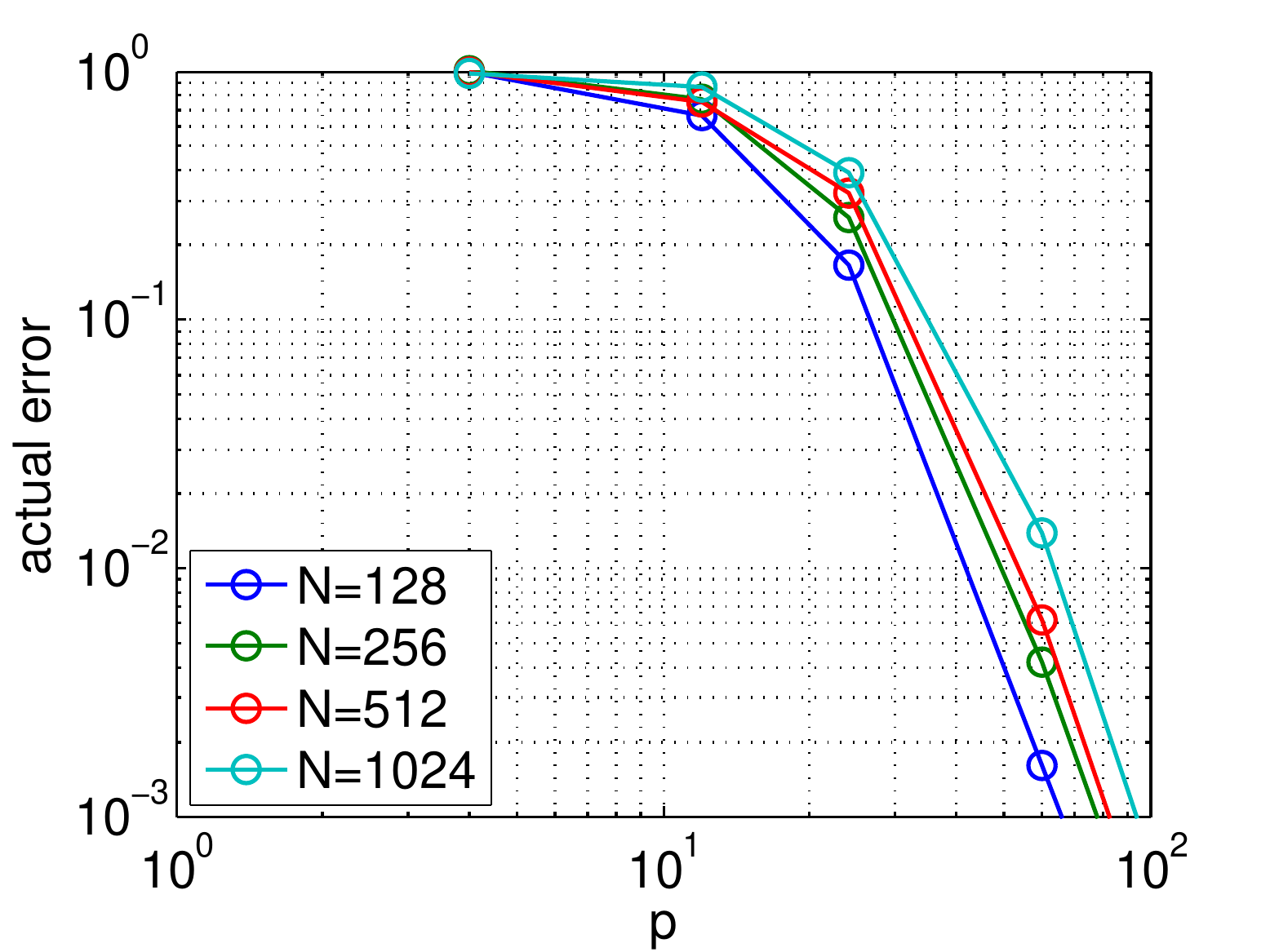}
\caption{Probing error of the $(1,1)$ block of the DtN map for the slow disk, fixed FD error level of $10^{-1}$, increasing $N$. This is the worst case, where $p$ follows a weak power law with $N$.}
\label{fig:c5pvsn}
\end{minipage}
\hspace{1cm}
\begin{minipage}[t]{0.45\linewidth}
\includegraphics[scale=.45]{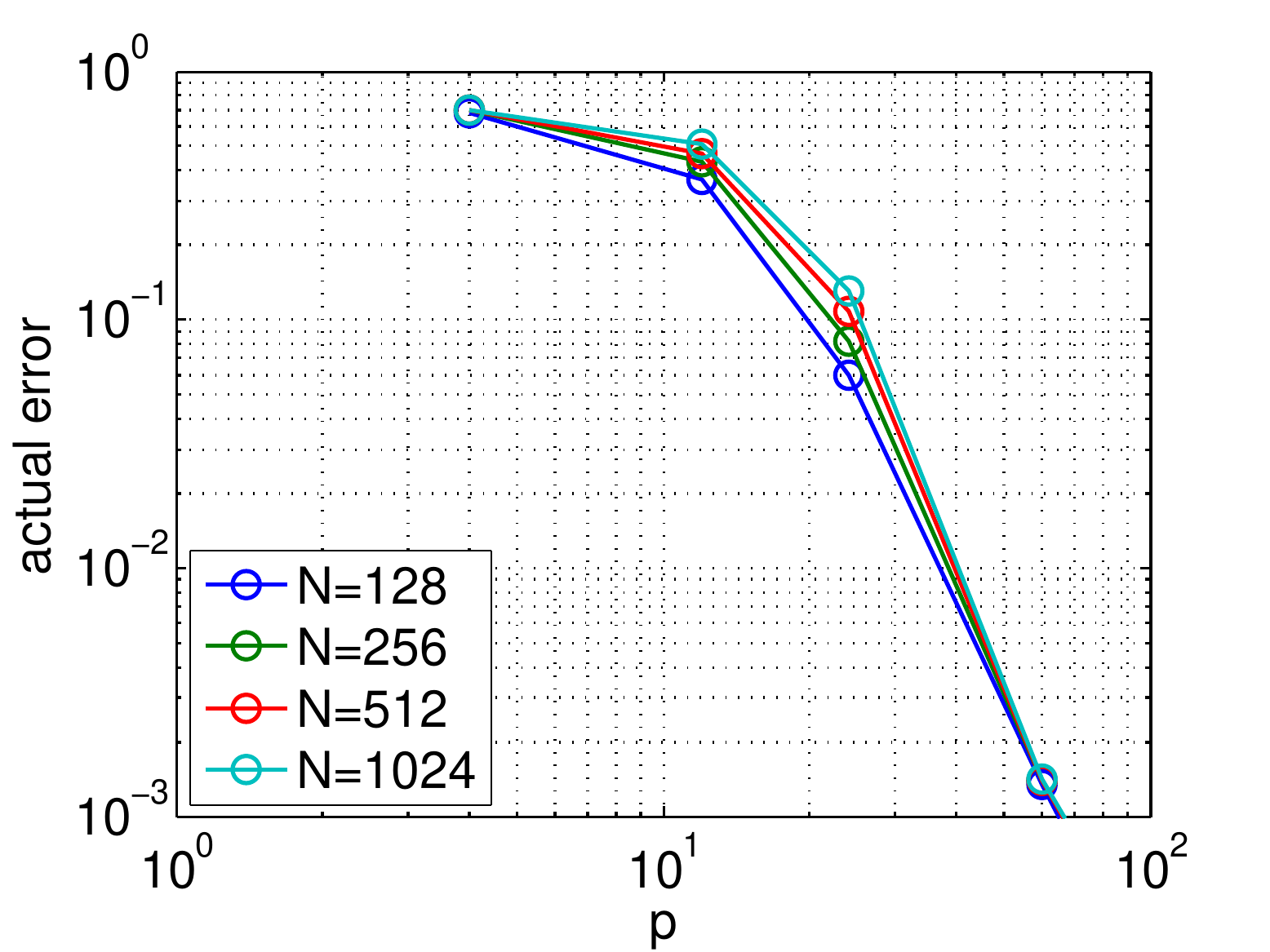}
\caption{Probing error of the $(1,1)$ block of the DtN map for the waveguide, fixed FD error level of $10^{-1}$, increasing $N$. Here $p$ follows a logarithmic law.}
\label{fig:c3pvsn}
\end{minipage}
\end{figure}

We now discuss how the number of basis matrices $p$ needed to achieve a desired accuracy depends on $N$ or $\omega$. To do this, we pick 4 consecutive powers of 2 as values for $N$, and find the appropriate $\omega$ such that the finite discretization error remains constant at $10^{-1}$, so that in fact $N \sim \omega^{1.5}$ as we have previously mentioned. We then probe the $(1,1)$ block of the corresponding DtN map, using the same parameters for all $N$, and observe the required $p$ to obtain a fixed probing error. The worst case we have seen in our experiments came from the slow disk. As we can see in Figure \ref{fig:c5pvsn}, $p$ seems to follow a very weak power law with $N$, close to $p \sim 15N^{.12}$ for a probing error of $10^{-1}$ or $p \sim 15N^{.2}$ for an probing error of $10^{-2}$. In all other cases, $p$ is approximately constant with increasing $N$, or seems to follow a logarithmic law with $N$ as for the waveguide (see Figure \ref{fig:c3pvsn}).

\section{Convergence of basis matrices for probing the half-space DtN map}\label{sec:pf}

In a uniform medium $c=1$, we have $K(r) = \frac{ik}{2r} H_1^{(1)}(kr)$. Its counterpart without the leading oscillatory factor is
\beq\label{eq:Hofr}
H(r) = \frac{ik}{2r} H_1^{(1)}(kr) e^{-ikr}.
\eeq
The expression of $K(r)$ is obtained by taking the mixed normal derivative $\pd_{n_x} \pd_{n_y} G(x,y)$ of the half-space Green's function, as explained for instance in \cite{FollandIntroPDE}, p. 92. (Green's representation formula involves $\pd_{n_y}G(x,y)$ to map Dirichlet data to the solution inside the domain, then the DtN map is asking to take an additional $\pd_{n_x}$.)

\subsection{Proof of Theorem \ref{teo:main}}

The domain of interest for the $r$ variable is $[r_0,1]$. Since the proof details will not depend on the constant $C$ in $r_0 = \frac{C}{k}$, we let $r_0 = \frac{1}{k}$ without loss of generality.

Expanding $K(r)$ in the system of Theorem \ref{teo:main} is equivalent to expanding $H(r)$ in polynomials of $r^{-1/\alpha}$ over $[r_0,1]$. it will be useful to perform the affine rescaling 
\[
\xi(r) = \frac{2}{r_0^{-1/\alpha} - 1} (r^{-1/\alpha} - 1 ) - 1 \qquad \Leftrightarrow \qquad r(\xi) = \left( \frac{\xi+1}{2} (r_0^{-1/\alpha} - 1) + 1 \right)^{-\alpha}
\]
so that the bounds $r \in [r_0,1]$ turn into $\xi \in [-1,1]$. We further write $\xi = \cos \theta$ with $\theta \in [0,\pi]$. Our strategy is to expand $H$ in Chebyshev polynomials $T_n(\xi)$. By definition, the best $p$-term approximation of $H(r)$ in polynomials of $r^{-1/\alpha}$ (best in a uniform sense over $[r_0,1]$) will result in a lower uniform approximation error than that associated with the $p$-term approximation of $H(r(\xi))$ in the $T_n(\xi)$ system. Hence in the sequel we overload notations and let $H_p$ for the $p$-term approximant of $H$ in our Chebyshev system.

We write out the Chebyshev series for $H(r(\xi))$ as 
$$H(r(\xi)) = \sum^{\infty}_{j=0} c_j T_j(\xi), \ c_j = \frac{2}{\pi} \int_{-1}^1 \frac{H(r(\xi)) T_j(\xi)}{(1-\xi^2)^{1/2}} \ d\xi, $$
with $T_j(\xi)=\cos{(j(\cos^{-1}\xi))}$, and $c_j$ alternatively written as
$$ c_j = \frac{2}{\pi} \int_0^\pi H(r(\cos{\theta})) \cos{j\theta} \ d \theta = \frac{1}{\pi} \int_0^{2\pi} H(r(\cos{\theta})) \cos{j\theta} \ d \theta. $$
The expansion will converge fast because we can integrate by parts in $\theta$ and afford to take a few derivatives of $H$, say $M$ of them, as done in \cite{tadmor}. After noting that the boundary terms cancel out because of periodicity in $\theta$, we express the coefficients $c_j$ for $j > 0$, up to a sign, as
$$ c_j  = \pm \frac{1}{\pi j^M} \int_0^{2\pi} \sin{j\theta} \frac{d^M}{d\theta^M} H(r(\cos{\theta})) \ d\theta, \qquad \ M \ \text{odd,} $$
$$ c_j  = \pm \frac{1}{\pi j^M} \int_0^{2\pi} \cos{j\theta} \frac{d^M}{d\theta^M} H(r(\cos{\theta})) \ d\theta, \qquad \ M \ \text{even.} $$
It follows that, for $j > 0$, and for all $M > 0$,
\[
\left| c_j \right| \leq \frac{2}{j^M} \max_\theta \left| \frac{d^M}{d\theta^M} H(r(\cos{\theta})) \right|.
\]
Let $B_M$ be a bound on this $M$-th order derivative. The uniform error we make by truncating the Chebyshev series to $H_p=\sum^{p}_{j=0} c_j T_j$ is then bounded by %
$$ \left\| H-H_p \right\|_{L^\infty[r_0,1]} \leq \sum_{j=p+1}^\infty \left|c_j \right| \leq 2 B_M \sum_{j=p+1}^\infty \frac{1}{j^{M}} \leq \frac{2B_M}{(M-1) p^{M-1}}, \qquad \ p > 1. $$
The final step is a simple integral comparison test.

The question is now to find a favorable estimate for $B_M$, from studying successive $\theta$ derivatives of $H(r)$ in (\ref{eq:Hofr}).  One of us already derived a bound for the derivatives of Hankel functions in \cite{flatland}: given any $C > 0$, we have
\begin{equation}\label{eq:derivH}
\left| \frac{d^m}{dr^m} \left( H_1^{(1)}(kr) e^{-ikr} \right) \right| \leq C_m (kr)^{-1/2} r^{-m} \qquad \ \text{for} \ kr \geq C.
\end{equation}
The change of variables from $r$ to $\theta$ results in
\begin{eqnarray*}
\frac{dr}{d\theta}&=& \frac{d\xi}{d\theta} \frac{dr}{d\xi} = \left(-\sin \theta \right) \left( -\alpha  \left( \frac{\xi+1}{2} (r_0^{-1/\alpha} - 1) + 1 \right)^{-\alpha-1} \frac{\left(r_0^{-1/\alpha}-1\right)}{2} \right) \\
&=&\left(-\sin \theta \right) \left( -\alpha \ r^{1+1/\alpha} \ \frac{k^{1/\alpha}(1-r_0^{1/\alpha})}{2} \right) = r (kr)^{1/\alpha} \ \frac{\alpha \sin \theta (1-r_0^{1/\alpha})}{2}.
\end{eqnarray*}
Derivatives of higher powers of $r$ are handled by the chain rule, resulting in
\[
\frac{d}{d \theta} (r^p) = p r^p (kr)^{1/\alpha} \ \frac{\alpha \sin \theta (1-r_0^{1/\alpha})}{2}.
\]
We see that the action of a $\theta$ derivative is essentially equivalent to multiplication by $(kr)^{1/\alpha}$. As for higher derivatives of powers of $r$, it is easy to see by induction that the product rule has them either hit a power of $r$, or a trigonometric polynomial of $\theta$, resulting in a growth of at most $(kr)^{1/\alpha}$ for each derivative:
\[
| \frac{d^m}{d\theta^m} r^p | \leq C_{m,p,\alpha} \, r^p (kr)^{m/\alpha}.
\]
These estimates can now be combined to bound $\frac{d^m}{d \theta^m} \left( H_1^{(1)}(kr) e^{-ikr} \right)$. One of two scenarios occur when applying the product rule:
\bit
\item either $\frac{d}{d\theta}$ hits $\frac{d^{m_2}}{d\theta^{m_2}} \left( H_1^{(1)}(kr) e^{-ikr} \right)$ for some $m_2 < m$. In this case, one negative power of $r$ results from $\frac{d}{dr}$, and a factor $r (kr)^{1/\alpha}$ results from $\frac{dr}{d\theta}$;
\item or $\frac{d}{d\theta}$ hits some power of $r$, or some $\frac{d^{m_1} r}{d\theta^{m_1}}$ for some $m_1 < m$, resulting in a growth of an additional factor $(kr)^{1/\alpha}$.
\eit
Thus, we get a $(kr)^{1/\alpha}$ growth factor per derivative in every case. The situation is completely analogous when dealing with the slightly more complex expression $\frac{d^m}{d \theta^m} \frac{1}{r} \left( H_1^{(1)}(kr) e^{-ikr} \right)$. The number of terms is itself at most factorial in $m$, hence we get
\[
| \frac{d^m}{d \theta^m} \frac{k}{r} \left( H_1^{(1)}(kr) e^{-ikr} \right) | \leq C_{m, \alpha} \frac{k}{r} (kr)^{\frac{m}{\alpha} - \frac{1}{2}} \leq C_{m,\alpha} k^2 (kr)^{\frac{m}{\alpha} - \frac{3}{2}}.
\]

We now pick $m \leq M = \lfloor 3 \alpha /2 \rfloor$, so that the max over $\theta$ is realized when $r = 1/k$, and $B_M$ is on the order of $k^2$. It follows that
$$ \left\| H-H_p \right\|_{L^\infty[r_0,1]} \leq  C_\alpha \ \frac{k^2}{p^{\lfloor 3\alpha/2 \rfloor - 1}}, \qquad \ p > 1, \ \alpha > 2/3.$$
The kernel of interest, $K(r) = H(r) e^{ikr}$ obeys the same estimate if we let $K_p$ be $p$-term approximation of $K$ in the Chebyshev system modulated by $e^{ikr}$.

We now turn to the operator norm of $\tilde{D} - \tilde{D}_p$ with kernel $\tilde{K} -\tilde{K}_p$, where $\tilde{K}(r) = K(r) \chi_{[r_0,1]}(r)$. Namely,
$$(\tilde{D}-\tilde{D}_p)g(x)=\int_0^1 (\tilde{K}-\tilde{K}_p)(|x-y|)g(y) \ dy. $$
We use the Cauchy-Schwarz inequality to bound
\begin{eqnarray*}
\|(\tilde{D}-\tilde{D}_p)g\|_2 &=& \left(\int_{0\leq x \leq 1}  \left|\int_{0\leq y \leq 1, \ |x-y|\geq r_0} (K-K_p)(|x-y|)g(y) \ dy \right|^2 dx \right)^{1/2} \\
& \leq & \left(\int_{0\leq x \leq 1} \int_{0\leq y \leq 1, \ |x-y|\geq r_0}\left|(K-K_p)(|x-y|)\right|^2 \ dy dx \right)^{1/2} \|g\|_2 \\
& \leq & \left( \int_{0\leq x \leq 1} \int_{0\leq y \leq 1, \ |x-y|\geq r_0} 1 \ dy \ dx \right)^{1/2} \|g\|_2 \ \max_{0\leq x,y \leq 1, \ |x-y| \geq r_0} |(K-K_p)(|x-y|)|  \\
& \leq &  \|g\|_2 \ \| K-K_p \|_{L^{\infty}[r_0,1]}.
\end{eqnarray*}
Assembling the bounds, we have
\[
\| \tilde{D}-\tilde{D}_p \| \leq C_{\alpha} \, p^{1 - \lfloor 3 \alpha / 2 \rfloor} \, k^2.
\]
It suffices therefore to show that $\| \tilde{K} \|_\infty = \| K \|_{L^\infty[r_0,1]}$ is on the order of $k^2$ to complete the proof. Letting $z = kr$, we see that
$$\max_{r_0 \leq r \leq 1} |K(r)|=\frac{k^2}{2} \max_{1 \leq z \leq k} \left|\frac{H_1^{(1)}(z) }{z} \right| \geq C k^2.$$
The last inequality follows from the fact that there exist positive constants $c_1$, $c_2$ such that $c_1 z^{-3/2} \leq \left|H_1^{(1)}(z)/z \right|  \leq c_2 z^{-3/2}$, which one of us proved in
\cite{flatland}.

\subsection{Numerical confirmation}

\begin{figure}[ht]
\begin{minipage}[t]{0.48\linewidth}
\includegraphics[scale=.5]{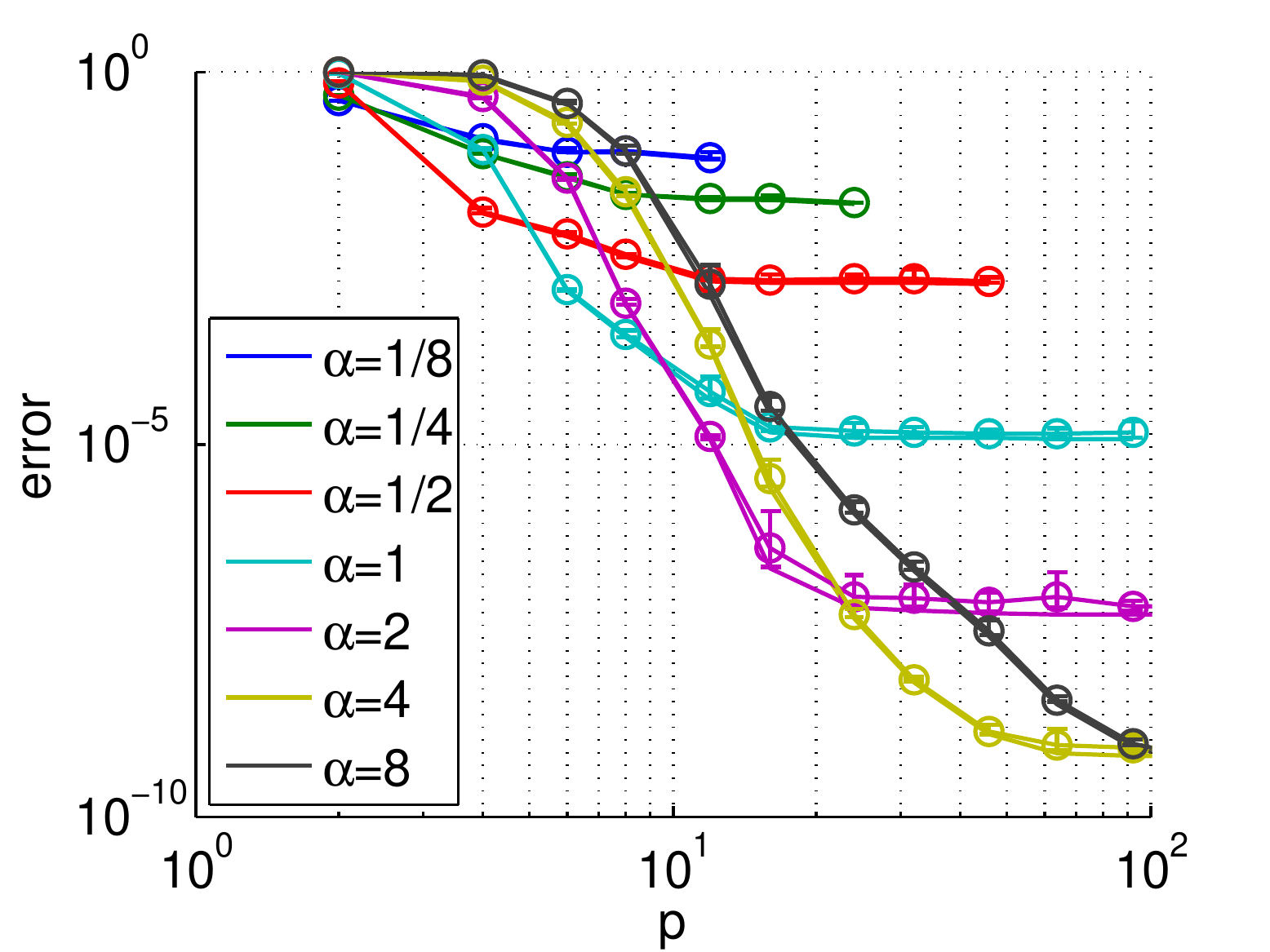}
\caption{Probing error of the half-space DtN map($q=1$, 10 trials, circle markers and error bars) compared to the approximation error (line), $c\equiv 1$, $L=1/4$, $\alpha=2$, $n=1024$, $\omega=51.2$.}
\label{q1}
\end{minipage}
\hspace{0.1cm}
\begin{minipage}[t]{0.48\linewidth}
\includegraphics[scale=.5]{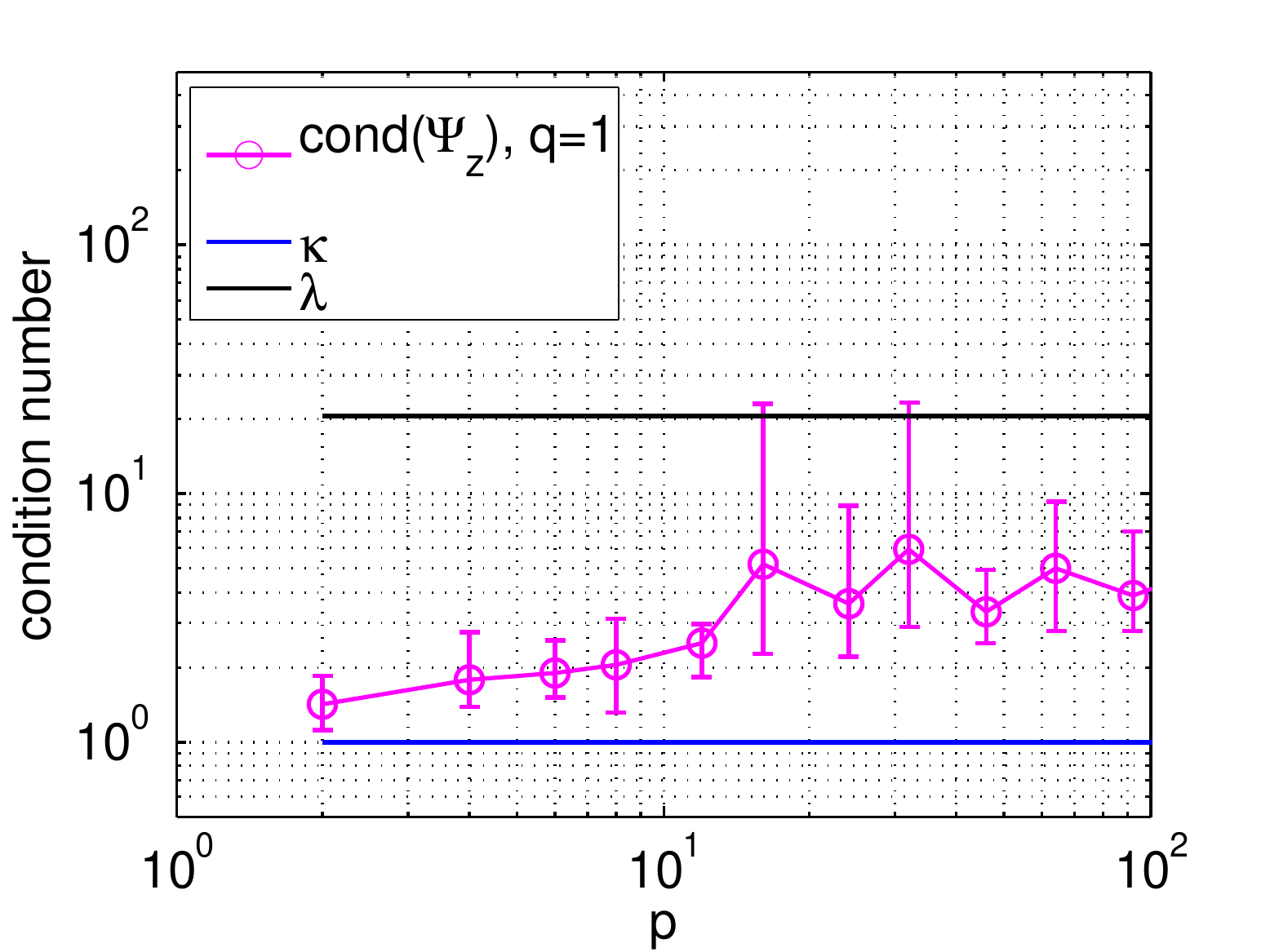}
\caption{Condition numbers for probing the half-space DtN map, $c\equiv 1$, $L=1/4$, $\alpha=2$, $n=1024$, $\omega=51.2$, $q=1$, 10 trials.}
\label{p2}
\end{minipage}
\end{figure}

In order to use Theorem \ref{teo:main} to obtain convergent basis matrices, we orthonormalize the set $\left\{(r)^{-j/\alpha}  \right\}_{j=0}^{p-1}$, put in oscillations, and use it as a basis for probing the DtN map.
Thus we obtain pre-basis matrices
$$(\beta_j)_{\ell m}=\frac{e^{ikh|\ell-m|}}{|\ell-m|^{j/\alpha}} \ \text{for} \ \ell \neq m,$$
with $(\beta_j)_{\ell \ell}=0$. We add to this set the identity matrix in order to capture the diagonal of $D$, and orthogonalize the resulting collection to get the $B_j$. Alternatively, we have noticed that orthogonalizing the set of $\beta_j$'s with
\begin{equation}\label{halfbasis}
 (\beta_j)_{\ell m}=\frac{e^{ikh|\ell-m|}}{(h+|\ell-m|)^{j/\alpha}}
\end{equation}
works just as well, and is simpler because there is no need to treat the diagonal separately. This is the same as we do for the exterior problem.

The convergent basis matrices in \eqref{halfbasis} have been used to obtain a numerical confirmation of Theorem \ref{teo:main}, again for the half-space DtN map. To obtain the DtN map in this setup, instead of solving the exterior problem with a PML on all sides, we solve a problem on a thin strip, with a random Dirichlet boundary condition (for probing) one of the long edges, and a PML on the other three sides. In Figure \ref{q1}, we show the approximation error, which we expect will behave as the approximation error. We also plot error bars for the probing error, corresponding to ten trials of probing, with $q=1$. The probing results are about as good as the approximation error, because the relevant condition numbers are all well-behaved as we see in Figure \ref{p2} for the value of choice of $\alpha=2$. Back to the approximation error, we notice in Figure \ref{q1} that increasing $\alpha$ delays the onset of convergence as expected, because of the factor which is factorial in $\alpha$ in the approximation error of Theorem \ref{teo:main}. And we can see that, for small $\alpha$, we are taking very high inverse powers of $r$, an ill-conditioned operation. Hence the appearance of a convergence plateau for smaller $\alpha$ is explained by ill-conditioning of the basis matrices, and the absence of data points is because of computational overflow.

Finally, increasing $\alpha$ from $1/8$ to $2$ gives a higher rate of convergence, as it should because in the error we have the factor $p^{-3\alpha/2}$, which gives a rate of convergence of $3\alpha/2$. This is roughly what we obtain numerically. As discussed, further increasing $\alpha$ is not necessarily advantageous since the constant $C_\alpha$ in \ref{teo:main} grows fast in $\alpha$.

\section{Discussion}

Probing the DtN map $D$ ultimately makes sense in conjunction with a fast algorithm for its application. In full matrix form, $D$ costs $N^2$ operations to apply. With the help of a compressed representation, this count becomes $p$ times the application complexity of any atomic basis function $B_j$, which may or may not be advantageous depending on the particular expansion scheme. The better solution for a fast algorithm, however, is to post-process the compressed expansion from probing into a slightly less compressed, but more algorithmically favorable one, such as h-matrix or butterfly. These types of matrix structure are not parametrized enough to lend themselves to efficient probing however -- see for instance \cite{Hmatvec} for an illustration of the large number of probing vectors required -- but will give rise to faster algorithms for matrix-vector multiplication. Hence the feasibilty of probing, and the availability of a fast algorithm for matrix-vector multiplcation, are two different goals that require different expansion schemes.

As for the complexity of solving the Helmholtz equation, reducing the PML confers the advantage of making the number of nonzeros in the matrix $L$ (of Section \ref{sec:strip}) independent of the width of the PML. After elimination of the layer, it is easy to see that $L$ has about $20N^2$ nonzero entries, instead of the $5N^2$ one would expect from a five-point stencil discretization of the Helmholtz equation, because the matrix $D$ (part of a small block of $L$) is in general full. Although obtaining a fast matrix-vector product for our approximation of $D$ could reduce the application cost of $L$ from $20N^2$ to something closer to $5N^2$, it should be noted that the asymptotic complexity does not change -- only the constant does, by a factor 4 at best. Hence the discussion about fast algorithms for $D$ is not as crucial as the idea of reducing the PML in the first place, when the goal is to apply $L$ fast in an iterative method.

\end{document}